# ASYMPTOTIC BEHAVIOR OF DIVERGENCES AND CAMERON–MARTIN THEOREM ON LOOP SPACES[1]

By Xiang Dong Li

*University of Oxford and Université Paul Sabatier*

*Dedicated to my grand father Li Xun-Cheng*

We first prove the $L^p$-convergence ($p \geq 1$) and a Fernique-type exponential integrability of divergence functionals for all Cameron–Martin vector fields with respect to the pinned Wiener measure on loop spaces over a compact Riemannian manifold. We then prove that the Driver flow is a smooth transform on path spaces in the sense of the Malliavin calculus and has an $\infty$-quasi-continuous modification which can be quasi-surely well defined on path spaces. This leads us to construct the Driver flow on loop spaces through the corresponding flow on path spaces. Combining these two results with the Cruzeiro lemma [*J. Funct. Anal.* **54** (1983) 206–227] we give an alternative proof of the quasi-invariance of the pinned Wiener measure under Driver's flow on loop spaces which was established earlier by Driver [*Trans. Amer. Math. Soc.* **342** (1994) 375–394] and Enchev and Stroock [*Adv. Math.* **119** (1996) 127–154] by Doob's $h$-processes approach together with the short time estimates of the gradient and the Hessian of the logarithmic heat kernel on compact Riemannian manifolds. We also establish the $L^p$-convergence ($p \geq 1$) and a Fernique-type exponential integrability theorem for the stochastic anti-development of pinned Brownian motions on compact Riemannian manifold with an explicit exponential exponent. Our results generalize and sharpen some earlier results due to Gross [*J. Funct. Anal.* **102** (1991) 268–313] and Hsu [*Math. Ann.* **309** (1997) 331–339]. Our method does not need any heat kernel estimate and is based on quasi-sure analysis and Sobolev estimates on path spaces.

**Contents**

1. Introduction

Received August 2002; revised October 2003.
[1]Supported in part by a Post-doctoral Research Fellowship of the University of Oxford.

2. The Airault–Malliavin–Sugita–Watanabe inequality
3. Sobolev estimates of divergence functionals on path spaces
4. The $L^p(\nu)$-convergence of divergence functionals
5. Exponential integrability of divergence functionals
6. Smoothness of Driver's flow on path spaces
7. Cameron-Martin theorem on loop spaces
8. Stochastic anti-development of pinned Brownian motions
9. Two remarks on Doob's $h$-processes approach

**1. Introduction.** Let $M$ be a compact connected Riemannian manifold and let $m_0 \in M$ be a fixed point. Let $\nabla$ be the Levi–Civita connection on $M$. The orthonormal frame bundle over $M$ is denoted by $O(M)$. The path space (resp., the loop space) over $M$ is defined by $P_{m_0}(M) = \{\gamma \in C([0,1], M) : \gamma(0) = m_0\}$ (resp., $L_{m_0}(M) = \{\gamma \in P_{m_0}(M) : \gamma(1) = m_0\}$). Let $H$ be the $\mathbb{R}^d$-valued Cameron–Martin space, that is, the set of absolutely continuous functions $h : [0,1] \to \mathbb{R}^d$ such that $h(0) = 0$ and $\dot{h} \in L^2([0,1])$. Let $H_0$ be its subspace with zero values at time 1, that is, $H_0 = \{h \in H : h(1) = 0\}$.

Let $\Delta$ be the Laplace–Beltrami operator on $M$. The Wiener measure on $P_{m_0}(M)$, denoted by $\mu$, is the law of $M$-valued Brownian motion (with generator $\Delta/2$) starting at $m_0$. The pinned Wiener measure on $L_{m_0}(M)$, denoted by $\nu$, is the law of the conditional Brownian motion $\{\gamma(s), s \in [0,1]\}$ on $M$ such that $\gamma(0) = \gamma(1) = m_0$. Intuitively, we have

$$\nu(\cdot) = \mu(\cdot|\gamma(1) = m_0).$$

Rigorously, if $p_t(x, y)$ denotes the heat kernel on $M$, then for any $\alpha < 1$,

$$(1.1) \qquad \frac{d\nu}{d\mu}\bigg|_{\mathcal{F}_\alpha}(\gamma) = \frac{p_{1-\alpha}(\gamma(\alpha), m_0)}{p_1(m_0, m_0)},$$

where $\mathcal{F}_\alpha = \sigma(\gamma(s), s \in [0, \alpha])$. For details, see, for example, [2] and [7].

For $\mu$-a.s. $\gamma \in P_{m_0}(M)$ [resp., $\nu$-a.s. $\gamma \in L_{m_0}(M)$], one can use the Itô SDE theory to define the stochastic parallel transport $U_s(\gamma) : T_{m_0}M \to T_{\gamma(s)}M$ as the unique $O(M)$-valued stochastic process satisfying the following covariant SDE:

$$\nabla_{\circ d\gamma(s)} U_s(\gamma) = 0,$$

with the initial condition $U_0(\gamma) = Id_{T_{m_0}M}$, where $Id_{T_{m_0}M}$ is the identity transform over $T_{m_0}M$. See, for example, [2]. For all $h \in H$ (resp., $h \in H_0$), the Cameron–Martin vector field $D_h$ on $P_{m_0}(M)$ [resp., $L_{m_0}(M)$] is defined by: for $\mu$-a.s. $\gamma \in P_{m_0}(M)$ [resp., $\nu$-a.s. $\gamma \in L_{m_0}(M)$],

$$(1.2) \qquad D_h(\gamma)(s) = U_s(\gamma)h(s) \quad \forall s \in [0,1].$$

In [6, 8, 14], the classical Cameron–Martin theorem has been generalized to the path space $(P_{m_0}(M), \mu)$. That is to say, for all fixed $h \in H$, the vector field $D_h$ generates a global flow $\{\Phi_t, t \in \mathbb{R}\}$ which can be $\mu$-a.s. well defined on $P_{m_0}(M)$ (see Section 6) under which the Wiener measure is quasi-invariant and an integration by parts formula holds. Concerning the same issue on the loop space $(L_{m_0}(M), \nu)$, Driver [7] proved that for any $h \in \mathcal{C}^1 \cap H_0$ (the set of Lipschitz Cameron–Martin vectors $h \in H_0$), the vector field $D_h$ generates a global flow $\{\widetilde{\Phi}_t : t \in \mathbb{R}\}$ on $L_{m_0}(M)$ such that, for $\nu$-a.s. $\gamma \in L_{m_0}(M)$,

$$\dot{\widetilde{\Phi}}_t(\gamma) = D_h(\widetilde{\Phi}_t(\gamma)), \qquad \widetilde{\Phi}_0(\gamma) = \gamma,$$

and the Wiener measure $\nu$ is quasi-invariant under the flow $\widetilde{\Phi}_t$; that is, the measure $(\widetilde{\Phi}_t)_*\nu$ is equivalent to $\nu$. Moreover, an integration by parts formula holds on $(L_{m_0}(M), \nu)$: for two cylindrical functionals $F$ and $G$ on $L_{m_0}(M)$, we have

$$E_\nu(D_h F G) = E_\nu(F(-D_h G + \delta(h)G))$$

and

$$(1.3) \qquad \frac{d(\widetilde{\Phi}_t)_*\nu}{d\nu}(\gamma) = \exp\left(\int_0^t \delta(h)(\widetilde{\Phi}_{-s}(\gamma))\, ds\right),$$

where $\delta(h)$ is the so-called divergence functional on $(L_{m_0}(M), \nu)$ defined as follows: for $\nu$-a.s. $\gamma \in L_{m_0}(M)$,

$$(1.4) \qquad \delta(h)(\gamma) = \int_0^1 (\dot{h}(s) + \tfrac{1}{2}\mathrm{Ric}_{U_s(\gamma)}(h(s)), dx(s)).$$

Here Ric denotes the Ricci curvature form over $O(M)$ and $\{x(s), s \in [0,1]\}$ is the stochastic anti-development of $\{\gamma(s), s \in [0,1]\}$, denoted by $x = I^{-1}(\gamma)$, and is given by the following Stratonovich stochastic integral:

$$(1.5) \qquad x(s) = \int_0^s U_r^{-1}(\gamma) \circ d\gamma_r, \qquad s \in [0, 1].$$

The complete theory of integration by parts formula on the loop space $(L_{m_0}(M), \nu)$ for all the Cameron–Martin vector fields $D_h$ with $h \in H_0$ was first proved by Enchev and Stroock [9], where the authors have also proved the quasi-invariance of the pinned Wiener measure under the flow generated by $D_h$. In [15] and [16], Hsu gave another approach to integration by parts formula on loop spaces which avoids the problem of the quasi-invariance of the pinned Wiener measure on the loop space. Let us mention that all the approaches appearing in [7, 9, 15, 16] relied strongly on the short time upper bound estimates on the gradient and the Hessian of the logarithmic of the heat kernel, and all these authors used the Doob $h$-processes method for conditional Brownian motion on a compact Riemannian manifold.

The purpose of this paper is to study the asymptotic behavior of divergence functionals, Driver flow and Cameron–Martin theorem on loop spaces as well as some related problems by a different approach. First, we will use the Airault–Malliavin–Sugita–Watanabe inequality (see Section 2) and some Sobolev estimates on the divergence functionals on the path space (see Section 3) to prove the $L^p$-convergence ($p \geq 1$) and the Fernique-type exponential integrability of divergence functionals with respect to the pinned Wiener measure on loop spaces; see Sections 4 and 5. Second, we will prove that the Driver flow $\{\Phi_t, t \in \mathbb{R}\}$ is a smooth transform on path spaces in the sense of the Malliavin calculus and has $\infty$-quasi-continuous version denoted by $\{\widetilde{\Phi}_t, t \in \mathbb{R}\}$ which can be quasi-surely well defined up to a slim subset of the path space $P_{m_0}(M)$. Moreover, we prove that if $h \in H_0$, then $\{\widetilde{\Phi}_t, t \in \mathbb{R}\}$ actually realizes the Driver flow generated by the vector field $D_h$ on the loop space $(L_{m_0}(M), \nu)$; see Sections 6 and 7. Third, we will combine these two results with the Cruzeiro lemma [4, 24] to give an alternative approach to the complete Cameron–Martin theorem on the loop space $(L_{m_0}(M), \nu)$ avoiding use of any heat kernel estimate; see Section 7. Finally, we use our method to establish the $L^p(\nu)$-convergence ($p \geq 1$) and a Fernique-type exponential integrability theorem for the stochastic anti-development of pinned Brownian motions on a compact Riemannian manifold equipped with any torsion-skew symmetric (TSS) connection; see Section 8. Our results generalize and sharpen some earlier results due to [13, 15, 16]. Our method is inspired by [25], where the authors first established the quasi-invariance of the pinned Wiener measure on the loop group over a compact Lie group under the left or the right action of a finite energy loop (which is nothing else than the Driver flow on the loop group). In some sense, it leads us to get *sharp or better* estimate than the direct approach based on Doob's $h$-theory and heat kernel estimates: see Section 8 and Section 9.

To state our main results, let us follow Hsu [15, 16] to introduce a sequence of functionals as follows: for any $h \in H$ and $s < 1$, let

$$\delta_s(h)(\gamma) = \int_0^s (\dot{h}(s) + \tfrac{1}{2}\operatorname{Ric}_{U_s(\gamma)}(h(s)), dx(s)), \qquad x = I^{-1}(\gamma).$$

Note that $\delta_s(h)$ is $\mu$-a.s. well defined on $P_{m_0}(M)$. Since $\nu$ is equivalent to $\mu$ on $\mathcal{F}_s$, compare (1.1), we see that $\delta_s(h)$ is also well defined for $\nu$-a.s. $\gamma \in L_{m_0}(M)$.

Now we are in a position to state our main results of this paper.

THEOREM 1.1. *Let $h \in H_0$. Then the divergence functional $\delta(h)$ [formally given by (1.4)] can be realized as the $L^p(\nu)$-limit of $\delta_s(h)$ as $s \to 1$ for all $p \geq 1$. In fact, for all $p \geq 1$, there is a constant $C_p$ such that, for all $h \in H_0$,*

$$\|\delta_s(h) - \delta(h)\|_{L^p(\nu)} \leq C_p \left(\int_s^1 |\dot{h}(t)|^2\, dt\right)^{1/2}$$

*and*

$$\|\delta(h)\|_{L^p(\nu)} \leq C_p \|h\|_{H_0}.$$

*Moreover, for all*

$$\lambda < \lambda_0 = \frac{1}{(2 + \|\operatorname{Ric}\|_\infty)\|h\|_H},$$

*we have*

$$E_\nu[\exp(\lambda|\delta(h)|^2)] < +\infty,$$

*or, equivalently,*

$$\lim_{t \to \infty} \frac{1}{t^2} \log \nu(\{\gamma \in L_{m_0}(M) : |\delta(h)| > t\}) \leq -\frac{1}{(2 + \|\operatorname{Ric}\|_\infty)\|h\|_H}.$$

COROLLARY 1.2. *For all $p > 1$, the gradient operator $D$ on the loop space $L_{m_0}(M)$ is closable from $L^p(L_{m_0}(M), \nu)$ into $L^p(L_{m_0}(M), H_0, \nu)$.*

COROLLARY 1.3. *For all $h \in H_0$ and all $\lambda > 0$, we have*

$$E_\nu[\exp(\lambda|\delta(h)|)] < \infty.$$

THEOREM 1.4. *For any $h \in H$, there exists an $\infty$-quasi-continuous version of the Driver flow $\{\Phi_t, t \in \mathbb{R}\}$ which can be well defined up to a slim subset of $P_{m_0}(M)$.*

Let $\{\widetilde{\Phi}_t, t \in \mathbb{R}\}$ be a fixed $\infty$-quasi-continuous version of $\{\Phi_t, t \in \mathbb{R}\}$. By the disintegration principle of the Wiener measure, $\{\widetilde{\Phi}_t, t \in \mathbb{R}\}$ can be $\nu$-a.s. well defined on $L_{m_0}(M)$.

THEOREM 1.5. *Let $h \in H_0$. Then $\{\widetilde{\Phi}_t, t \in \mathbb{R}\}$ is the flow generated by the vector field $D_h$ on $L_{m_0}(M)$. Moreover, the pinned Wiener measure $\nu$ is quasi-invariant under $\{\widetilde{\Phi}_t, t \in \mathbb{R}\}$. More precisely, if we let*

$$K_t := \frac{d(\widetilde{\Phi}_t)_*\nu}{d\nu},$$

*then*

$$K_t = \exp\left(\int_0^t \delta(h)(\widetilde{\Phi}_{-s}(\gamma))\,ds\right), \qquad \nu\text{-a.s. } \gamma \in L_{m_0}(M),$$

*and for all $p > 1$ with the conjugate exponent $q$, that is, $\frac{1}{p} + \frac{1}{q} = 1$, we have*

$$\|K_t\|^p_{L^p(\nu)} \leq E_\nu[e^{pqt|\delta(h)|}].$$

REMARK 1.1. All the above results [as well as the Sobolev norms and capacities comparison inequalities (6.17) and (6.18) in Section 6] remain true if we replace the Levi–Civita connection by any torsion skew-symmetric (TSS) connection. In this case, we need only to replace the Ricci curvature Ric of the Levi–Civita connection by the Ricci curvature $\widehat{\mathrm{Ric}}$ of the dual connection $\widehat{\nabla}$ given by

$$\widehat{\nabla}_X Y = \nabla_X Y - T(X, Y), \qquad X, Y \in \Gamma(TM),$$

where $T$ is the torsion of our given TSS connection $\nabla$. Indeed, we have announced Theorem 1.1 in [20] (without giving the precise value $\lambda_0$) in this setting with an equivalent expression of the divergence functional $\delta(h)$ as used in [7]. In particular, we recapture the Malliavin–Malliavin theorem on the quasi-invariance of the pinned Wiener measure on loop group over a compact Lie group.

The following result generalizes and sharpens some earlier results due to [13], where the author established the $L^p(\nu)$-convergence and a Fernique-type exponential integrability theorem for the stochastic anti-development of pinned Brownian motions on a compact Lie group.

THEOREM 1.6. *Let $M$ be a compact Riemannian manifold equipped with a TSS connection, $m \in \mathbb{N}$, $m \geq 2$, $\alpha \in (\frac{1}{2m}, \frac{1}{2})$. Then for any*

$$\lambda < \lambda_0 := \tfrac{1}{2} \inf\{\|w\|_H^2 : w \in X, \ \|w\|_{2m,\alpha} = 1\},$$

*we have*

$$E_\nu[\exp(\lambda \|x\|_{2m,\alpha}^2)] < +\infty,$$

*where*

$$\|x\|_{2m,\alpha} = \left[\int_0^1 \int_0^1 \frac{\|x(t) - x(s)\|_{\mathbb{R}^d}}{|t-s|^{1+2m\alpha}} \, dt \, ds\right]^{1/2m}$$

*and*

$$x(s) = \int_0^1 U_s^{-1}(\gamma) \circ d\gamma(s), \qquad s \in [0,1], \ \nu\text{-a.s.} \ \gamma \in L_{m_0}(M).$$

*Moreover, for any $\lambda < \frac{1}{2}$, we have*

$$E_\nu\left[\exp\left(\lambda \sup_{s \in [0,1]} \|x(s)\|^2\right)\right] < +\infty.$$

*In addition, $x(s)$ converges to $x(1)$ in $L^p(\nu)$ for all $p \geq 1$ as $s$ tends to $1$ and there exists a constant $C_p$ such that*

$$\|x(s) - x(1)\|_{L^p(\nu)} \leq C_p (1-s)^{1/2}.$$

**2. The Airault–Malliavin–Sugita–Watanabe inequality.** Let $X = \{x \in C([0,1], \mathbb{R}^d) : x(0) = 0\}$ be the Wiener space, and let $\mu_0$ be the Wiener measure on $X$. For any $r \in \mathbb{N}$ and $p > 1$, we let $W^{r,p}(X)$ denote the $(r,p)$-Sobolev space on the Wiener space $X$ with the Sobolev norm $\|\cdot\|_{W^{r,p}(X)}$. Let $A_1, \ldots, A_d$ be the canonical horizontal vector fields on $O(M)$, and let $r_0 \in O(M)$ be a fixed orthonormal frame over $m_0$. Consider the horizontal SDE on $O(M)$:

$$dr_x(s) = \sum_{i=1}^d A_i(r_x(s)) \circ dx^i(s),$$

$$r_x(0) = r_0.$$

Let $\gamma_x(s) := \pi(r_x(s))$, $s \in [0,1]$. Then it is well known that $\{\gamma_x(s), s \in [0,1]\}$ is a Brownian motion on $M$ starting at $m_0$. The Wiener measure $\mu$ on the path space $P_{m_0}(M) = \{\gamma \in C([0,1], M) : \gamma(0) = m_0\}$ is given by the law of $\{\gamma_x(s), s \in [0,1]\}$, that is, $\mu = I_*\mu_0$, where $I : X \to P_{m_0}(M)$ is the Itô map given by

$$I(x) = \gamma_x, \qquad \mu_0\text{-a.s. } x \in X.$$

Consider the following $M$-valued Wiener functional $\Phi$:

$$\Phi(x) = \gamma_x(1).$$

By [1], $\Phi \in W^{\infty,\infty}(X, M)$ and $\Phi$ is nondegenerated, that is,

$$(\mathrm{Det}[\Phi(x)])^{-1} \in W^\infty(X, M),$$

where $W^{\infty,\infty}(X, M)$ is the set of all smooth $M$-valued Wiener functionals in the sense of Malliavin calculus, and

$$\mathrm{Det}[\Phi(x)] = \sqrt{\det[\nabla\Phi(x) \cdot \nabla\Phi(x)^\tau]},$$

where the determinant on the right-hand side is taken with respect to the Riemannian metric on $T_{\Phi(x)}M$ and $\nabla\Phi(x)^\tau$ denotes the adjoint of $\nabla\Phi(x) : H \to T_{\Phi(x)}M$.

Recall that if $f \in W^{\infty,\infty}(X, M)$, then $f$ has an $\infty$-quasi-continuous modification which can be well defined outside a slim subset of $X$. Moreover, if $f_1, f_2$ are two $\infty$-quasi-continuous modifications of $f$, then $f_1$ and $f_2$ only differ on a slim set. Let $\Phi^*$ be any quasi-continuous modification of $\Phi$. The following co-area formula is well known (see, e.g., [1, 23, 24]): there exists a family of area measures denoted by $\{da^y(\cdot), y \in M\}$ [where each $da^y(\cdot)$ is supported on the submanifold $S_y = \Phi^{*-1}(y)$] such that, for any $u \in W^\infty(X)$ and $v \in C^\infty(M)$,

$$\int_X u(x) v(\Phi(x)) [\mathrm{Det}\,\Phi](x) \, d\mu_0(x) = \int_M v(y) \int_{\Phi^{*-1}(y)} u^*(x) \, da^y(x) \, dy,$$

where $u^*$ denotes any $\infty$-quasi-continuous modification of $u$. Let

$$\nu_y(dx) = (\text{Det}[\Phi](x))^{-1} da^y(x).$$

Then for all $y \in M$, $\nu_y$ is a Borel probability measure supported on the submanifold $S_y = \Phi^{*-1}(y)$. Moreover, $\nu_y$ has no charge on any slim subset of $X$.

By [22, 33], the Itô map $I \colon X \to P_{m_0}(M)$ is smooth in the sense of the Malliavin calculus and has an $\infty$-quasi-continuous modification. Throughout this paper, we let $\widetilde{I}$ denote a fixed $\infty$-quasi-continuous modification of the Itô map $I$. Using the dyadic polygonal approximation of $M$-valued Brownian motion, and by a similar argument used in Section 4 in [30], one can prove that the stochastic anti-development map given by (1.5) has an $\infty$-quasi-continuous version (denoted by $\widetilde{I}^{-1}$) which can be quasi-surely well defined on $P_{m_0}(M)$ and satisfies $\widetilde{I} \circ I = Id_X$ quasi-surely (i.e., except on a slim set of $X$). By this and using the capacity comparison inequality due to the author [22], we can easily prove that, for two different versions of $\infty$-quasi-continuous modification of $I$, say $I_1$ and $I_2$, $I_1|_{S_{m_0}}(S_{m_0})$ only differs from $I_2|_{S_{m_0}}(S_{m_0})$ on a slim subset of $P_{m_0}(M)$. Indeed, if we let $S = \{x : I_1(x) \neq I_2(x)\}$, then $S$ is a slim set of $X$. Let $O \subset X$ be an open set containing $S$ and with capacity $C_{r,p}(O) < \varepsilon$ for all $r \in \mathbb{N}$ and $p > 1$. By the capacity comparison inequality between the path space and the Wiener space (see [22]), we have

$$\widehat{C}_{r,p}(I_i(O)) \leq \alpha C_{2r,p+1}(I_i^{-1} \circ I_i(O)), \qquad i = 1,2,$$

where $\alpha = \alpha(r,p)$ is a constant and $\widehat{C}_{r,p}$ is the $(r,p)$-capacity on the path space $P_{m_0}(M)$ (for its definition, see [22]). Note that $I_i^{-1} \circ I_i = Id_X$ holds quasi-surely on $X$, $i=1,2$. Hence $C_{2r,p+1}(I_i^{-1} \circ I_i(O)) = C_{2r,p+1}(O) < \varepsilon$, $i=1,2$. Since $\varepsilon$ is arbitrary, we get $\widehat{C}_{r,p}(I_i(S)) = 0$, $\forall r \in \mathbb{N}$, $p > 1$, $i=1,2$. Thus, $\widetilde{I}|_{S_{m_0}}$, the restriction of $\widetilde{I}$ on the submanifold $S_{m_0} = \{x \in X : \gamma_x(1) = m_0\}$, is $\nu_{m_0}$-a.s. well defined. Moreover, $\widetilde{I}|_{S_{m_0}} \colon (S_{m_0}, \nu_{m_0}) \to (L_{m_0}(M), \nu)$ is a measure-theoretic isomorphism. That is to say, $\widetilde{I}|_{S_{m_0}}(S_0)$ only differs from $L_{m_0}(M)$ on a slim set of the path space $P_{m_0}(M)$ and

$$\nu = (\widetilde{I}|_{S_{m_0}})_* \nu_{m_0}.$$

The following result is due to [1, 33].

THEOREM 2.1. *There exist a pair $(p,r) \in (1,+\infty) \times \mathbb{N}$ and a constant $C > 0$ such that, for any $f \in W^\infty(X, \mathbb{R}^+)$, we have*

$$\int_{S_{m_0}} f^*(x) \nu_{m_0}(dx) \leq C \|f\|_{W^{r,p}(X)},$$

*where $f^*$ is any $\infty$-quasi-continuous modification of $f$.*

In fact, using the Watanabe generalized distribution theory on Wiener space, we can even specify the constant $C$ and the value of the pair $(r,p)$ appearing in Theorem 2.1 as follows. To this end, using the Nash–Whitney embedding theorem, we assume that $M$ is isometrically embedded into $\mathbb{R}^l$ with $l \geq d$.

THEOREM 2.2 (Airault–Malliavin–Sugita–Watanabe inequality). *For all* $p > 1$, $k \in \mathbb{N}$ *and* $f \in W^{2[l/2]+2+2k,p}(X)$, *we have*

$$\left| \int_{S_{m_0}} f^*(x)\, d\nu_{m_0}(x) \right| \leq C \|\delta_{m_0} \circ \Phi\|_{-2[l/2]-2-2k, p/(p-1)} \|f\|_{2[l/2]+2+2k, p},$$

*where*

$$C = [p_1(m_0, m_0)]^{-1}.$$

PROOF. For any $y \in M \subset \mathbb{R}^l$, let $\delta_y$ be the Dirac delta function at point $y$. Then $\delta \in \mathcal{S}_{-2r}(\mathbb{R}^l)$ for all $r \geq [\frac{l}{2}] + 1$, where $\mathcal{S}_{-2r}(\mathbb{R}^l)$ is the topological dual of $\mathcal{S}_{2r}(\mathbb{R}^l)$ (the completion of the Schwartz space $\mathcal{S}(\mathbb{R}^l)$ of rapidly decreasing $C^\infty$-functions on $\mathbb{R}^l$ by the norm $\|\cdot\|_{2r}$ defined by $\|\phi\|_{2r} = \|(1 + |x|^2 - \Delta)^r \phi\|_\infty$). See, for example, the proof of Theorem 4.2 in [32] and [34], Remark 2.2. Let $F \in W^{\infty,\infty}(X, M)$ be a smooth nondegenerate Wiener functional. Then for all $k = 0, 1, \ldots$ and $p > 1$, the map $y \in M \to \delta_y(F) \in W^{-2[l/2]-2-2k, p}(X)$ is $2k$-times continuous differentiable. Hence for any $f \in D^{2[l/2]+2+2k,p}(X)$ we have

$$\left| \int_X f(x) \delta_y(F(x))\, d\mu_0(x) \right| \leq \|\delta_y \circ F\|_{-2[l/2]-2-2k, p/(p-1)} \|f\|_{2[l/2]} + 2 + 2k, p.$$

In particular, taking $F = \Phi$ and using the fact that

$$\nu_{m_0}(dx) = \frac{\delta_{m_0}(\Phi(x))}{\int_X \delta_{m_0}(\Phi(x))\, d\mu_0(x)}\, d\mu_0(x),$$

we deduce the Airault–Malliavin–Sugita–Watanabe inequality with the constant $C$ given by (cf. [34])

$$C = \left[ \int_X \delta_{m_0}(\Phi(x))\, d\mu_0(x) \right]^{-1} = [p_1(m_0, m_0)]^{-1}. \qquad \square$$

**3. Sobolev estimates of divergence functionals on path spaces.** Following [22], for any $r \in \mathbb{N}$ and $p > 1$, we let $D^{r,p}(P_{m_0}(M))$ denote the $(r,p)$-Sobolev space on $P_{m_0}(M)$ with the Sobolev norm $\|\cdot\|_{D^{r,p}(P_{m_0}(M))}$ defined by

$$\|F\|_{D^{r,p}(P_{m_0}(M))} = \sum_{k=0}^{r} \|\|D^k F\|_{H^{\otimes k}}\|_p.$$

For any fixed $t \in [0,1]$, regarding $x \to r_x(t)$ as an $O(M)$-valued Wiener functional, we have $r.(t) \in W^{\infty,\infty}(X, O(M))$. More precisely, for any $n \in \mathbb{N}$ and any $h, h_1, \ldots, h_n \in H$, the following $H$-directional derivatives exist:

$$\nabla_h r_x(t) := \left\{\frac{d}{d\varepsilon} r_{x+\varepsilon h}(t)\right\}\bigg|_{\varepsilon=0},$$

$$\vdots$$

$$\nabla_{h_1 \ldots h_n} r_x(t) := \left\{\frac{D^n}{\partial \varepsilon_1 \cdots \partial \varepsilon_n} r_{x+\sum_{i=1}^n \varepsilon_i h_i}(t)\right\}\bigg|_{\varepsilon_1=\cdots=\varepsilon_n=0},$$

where $\frac{D}{\partial \varepsilon_i}$ denotes the Levi–Civita covariant derivative along the smooth curve $\varepsilon_i \longmapsto r_{x+\varepsilon_i h_i}(t)$ on $O(M)$. Moreover, we have the following proposition:

PROPOSITION 3.1 [22]. *There exist $D^{j_1,\ldots,j_n}_{s_1,\ldots,s_n} r_x(t) \in L^2([0,1]^n, T_{r_x(t)}O(M))$ such that:*

(i) $D^{j_1,\ldots,j_n}_{s_1,\ldots,s_n} r_x(t)$ *is adapted with respect to $\mathcal{F}_t = \sigma(x(s), s \in [0,t])$ for any fixed $s_1, \ldots, s_n \in [0,1]$ and $0 \leq j_1, \ldots, j_n \leq n$. Moreover,*

$$D^{j_1,\ldots,j_n}_{s_1,\ldots,s_n} r_x(t) = 0 \quad \text{if } s_1 \vee \cdots \vee s_n \in [t,1];$$

(ii) *for any $h_1, \ldots, h_n \in H$,*

$$\langle \nabla^n r_x(t), h_1 \otimes \cdots \otimes h_n \rangle_{H^{\otimes n}}$$
$$= \sum_{j_1,\ldots,j_n} \int_0^t \cdots \int_0^t D^{j_1,\ldots,j_n}_{s_1,\ldots,s_n} r_x(t) \dot{h}_1^{j_1}(s_1) \cdots \dot{h}_n^{j_n}(s_n) \, ds_1 \cdots ds_n;$$

(iii) *for any $p \geq 1$, we have*

$$(3.1) \qquad \sup_{s_1,\ldots,s_n \in [0,1]} E\left[\sup_{s_1 \vee \cdots \vee s_n \leq s \leq 1} \|D^{j_1,\ldots,j_n}_{s_1,\ldots,s_n} r_x(s)\|^p\right] < +\infty,$$

*where $\|D^{j_1,\ldots,j_n}_{s_1,\ldots,s_n} r_x(s)\|$ denotes the Riemannian norm of the vector field $D^{j_1,\ldots,j_n}_{s_1,\ldots,s_n} r_x(s)$ with respect to the Sasake Riemannian metric on $O(M)$ (for its definition, see [22]).*

Let $e_1, \ldots, e_d$ be the standard orthonormal basis of $\mathbb{R}^d$, and let $R$ be the Riemannian curvature tensor on $M$. For any $r \in O(M)$, the Ricci curvature over the frame $r$ is a real matrix given by

$$\mathrm{Ric}_r(a) = \sum_{i=1}^d r^{-1} \circ R(re_i, ra) \circ re_i \qquad \forall a \in \mathbb{R}^d.$$

Let

$$J(x,t) = \tfrac{1}{2} \mathrm{Ric}_{r_x(t)}.$$

By the chain rule and Proposition 3.1, we have $J(x,t) \in W^\infty(X, M(d,d))$, where $M(d,d)$ denotes the set of all $d \times d$ real matrices. Moreover, we have the following proposition:

PROPOSITION 3.2. *The Malliavin derivatives $D^{j_1,\ldots,j_n}_{s_1,\ldots,s_n} J(x,t)$ belong to $L^2([0,1]^n, M(d,d))$ and are adapted with respect to $\mathcal{F}_t = \sigma(x(s), s \in [0,t])$ for any fixed $s_1,\ldots,s_n \in [0,1]$ and any $0 \leq j_1,\ldots,j_n \leq n$. For any $h_1,\ldots,h_n \in H$,*

$$\langle \nabla^n J(x,t), h_1 \otimes \cdots \otimes h_n \rangle$$
$$= \int_0^t \cdots \int_0^t D^{j_1,\ldots,j_n}_{s_1,\ldots,s_n} J(x,t) \dot{h}_1^{j_1}(s_1) \cdots \dot{h}_n^{j_n}(s_n) \, ds_1 \cdots ds_n.$$

*Moreover, for any $p \geq 1$,*

$$(3.2) \qquad \sup_{s_1,\ldots,s_n \in [0,1]} E\left[ \sup_{s_1 \vee \cdots \vee s_n \leq t \leq 1} \|D^{j_1,\ldots,j_n}_{s_1,\ldots,s_n} J(x,t)\|^p_{H \cdot S} \right] < +\infty.$$

PROOF. The proof can be easily given by the chain rule and using Proposition 3.1. In particular, (3.1) yields (3.2). □

Let $h \in H$ and $D_h$ be the vector fields on $P_{m_0}(M)$ defined by (1.2). By integration by parts formula (see, e.g., [2, 6, 8, 10, 14]), the $L^2(\mu)$ adjoint of $D_h$ is given by $D_h^* = -D_h + \delta(h)$, where $\delta(h)$ is the divergence functional. Moreover, for $\mu$-a.s. $\gamma \in P_{m_0}(M)$, we have

$$\delta(h)(\gamma) = \int_0^1 (\dot{h}(\tau) + \tfrac{1}{2} \operatorname{Ric}_{r_x(\tau)} h(\tau), dx(\tau)), \qquad x = I^{-1}(\gamma).$$

Now we state the main result in this section.

THEOREM 3.3. *For any $r \in \mathbb{N}, p > 1$, there is a constant $C > 0$ such that, for all $h \in H$, we have*

$$\|\delta(h)\|_{D^{r,p}(P_{m_0}(M))} \leq C \|h\|_H.$$

PROOF. Let $\widehat{\delta(h)} = \delta(h) \circ I$. By the Sobolev norm comparison theorem (see [22]), we have

$$\|\delta(h)\|_{D^{r,p}(P_{m_0}(M))} \leq \alpha_{r,p} \|\widehat{\delta(h)}\|_{W^{2r,p+1}(X)},$$

where $\alpha_{r,p}$ is a constant. Hence we need only to prove that, for any $r \in \mathbb{N}, p \geq 2$,

$$(3.3) \qquad \|\widehat{\delta(h)}\|_{W^{r,p}(X)} \leq C \|h\|_H.$$

By induction and direct computation, it can be easily shown that, for any $n \geq 2$ and $h, h_1, \ldots, h_n \in H$, we have

$$D_s \widehat{\delta(h)}(x) = \dot{h}(s) + J(x,s)h(s) + \int_0^1 \langle D_s J(x,t)h(t), dx(t) \rangle,$$

$$(3.4) \quad D_{s_1,\ldots,s_n}^{j_1,\ldots,j_n} \widehat{\delta(h)}(x) = \int_0^1 \langle D_{s_1,\ldots,s_n}^{j_1,\ldots,j_n} J(x,t)h(t), dx(t) \rangle$$

$$+ \sum_{i=1}^n (D_{s_1,\ldots,\widehat{s_i},\ldots,s_n}^{j_1,\ldots,\widehat{j_i},\ldots,j_n} J(x,s_i)h(s_i))^{j_i},$$

where we use the notation $(a_1, \ldots, a_d)^j := a_j$, $j = 1, \ldots, n$.

By the Burkholder–Davis–Gundy inequality, we have $\|\widehat{\delta(h)}\|_p \leq C\|h\|_H$. It remains to prove that, for any $n \geq 1$, $p \geq 2$, there is a constant $C > 0$ such that

$$(3.5) \quad \|\|\nabla^n \widehat{\delta(h)}\|_{H \cdot S}\|_{L^p(\mu_0)} \leq C\|h\|_H.$$

Below we only give a proof of (3.5) for $n \geq 2, p \geq 2$. The proof for the case of $n = 1, p \geq 2$ is analogous. By definition, we have

$$\|\nabla^n \widehat{\delta(h)}\|_{H \cdot S}^p = \left[ \sum_{j_1,\ldots,j_n} \int_{[0,1]^n} |D_{s_1,\ldots,s_n}^{j_1,\ldots,j_n} \widehat{\delta(h)}|^2 \, ds_1 \cdots ds_n \right]^{p/2}$$

$$\leq 2^{p/2} \left[ \sum_{j_1,\ldots,j_n} \int_{[0,1]^n} \left| \int_0^1 \langle D_{s_1,\ldots,s_n}^{j_1,\ldots,j_n} J(x,t)h(t), dx(t) \rangle \right|^2 ds \right]^{p/2}$$

$$+ 2^{p/2} \left[ \sum_{j_1,\ldots,j_n} \sum_{i=1}^n \int_{[0,1]^n} |(D_{s_1,\ldots,\widehat{s_i},\ldots,s_n}^{j_1,\ldots,\widehat{j_i},\ldots,j_n} J(x,s_i)h(s_i))^{j_i}|^2 \, ds \right]^{p/2}$$

$$\leq I_1 + I_2.$$

By the Hölder inequality and the Burkholder–Davis–Gundy inequality, we have

$$I_1 = E \left[ \int_{[0,1]^n} \left| \int_0^1 \langle D_{s_1,\ldots,s_n}^{j_1,\ldots,j_n} J(x,t)h(t), dx(t) \rangle \right|^2 ds \right]^{p/2}$$

$$\leq \int_{[0,1]^n} E \left| \int_0^1 \langle D_{s_1,\ldots,s_n}^{j_1,\ldots,j_n} J(x,t)h(t), dx(t) \rangle \right|^p ds$$

$$\leq C \int_{[0,1]^n} E \left[ \int_0^1 |D_{s_1,\ldots,s_n}^{j_1,\ldots,j_n} J(x,t)h(t)|^2 \, dt \right]^{p/2} ds$$

$$\leq C\|h\|_H^p \int_{[0,1]^n} E \left( \sup_{t \in [\vee s_i, 1]} |D_{s_1,\ldots,s_n}^{j_1,\ldots,j_n} J(x,t)|^p \right) ds \qquad \text{[by (3.2)]}$$

$$\leq C\|h\|_H^p.$$

Similarly, we have

$$
\begin{aligned}
I_2 &= E\bigg[\int_{[0,1]^n}|(D^{j_1,\ldots,\widehat{j_i},\ldots,j_n}_{s_1,\ldots,\widehat{s_i},\ldots,s_n}J(x,s_i)h(s_i))^{j_i}|^2 ds\bigg]^{p/2}\\
&\leq E\bigg[\int_0^1\bigg(\sup_{s_i\in[\vee s_j,1]}|D^{j_1,\ldots,\widehat{j_i},\ldots,j_n}_{s_1,\ldots,\widehat{s_i},\ldots,s_n}J(x,s_i)|^2\bigg)h(s_i)\,ds_i\bigg]^{p/2}\\
&\leq \|h\|_H^p \sup_{s_j\in[0,1]} E\bigg[\sup_{s_i\in[\vee s_j,1]}|D^{j_1,\ldots,\widehat{j_i},\ldots,j_n}_{s_1,\ldots,\widehat{s_i},\ldots,s_n}J(x,s_i)|^p\bigg] \quad \text{[by (3.2)]}\\
&\leq C\|h\|_H^p.
\end{aligned}
$$

Combining the above inequalities for $I_1$ and $I_2$, we obtain (3.5) and hence (3.3). □

PROPOSITION 3.4. *For any $n\in\mathbb{N}, p>1$, there is a constant $C>0$ such that, for all $h\in H$, we have*

(3.6)
$$
\bigg\|\int_{s_1}^{s_2}\langle \dot h(t)+\tfrac{1}{2}\mathrm{Ric}_{r_x(t)}(h(t)), dx(t)\rangle\bigg\|_{D^{r,p}(P_{m_0}(M))}
$$
$$
\leq C\bigg[\bigg(\int_{s_1}^{s_2}|\dot h(t)|^2\,dt\bigg)^{1/2}+\bigg(\int_{s_1}^{s_2}|h(t)|^2\,dt\bigg)^{1/2}\bigg].
$$

PROOF. The proof is similar to the one of Theorem 3.3. □

COROLLARY 3.5. *The functional $\delta(h)$ has an $\infty$-quasi-continuous modification which can be $\nu$-a.s. well defined on $L_{m_0}(M)$.*

PROOF. By Theorem 3.3, $\delta(h)\in D^{\infty,\infty}(P_{m_0}(M))=\bigcap_{r\in N, p>1}D^{r,p}(P_{m_0}(M))$. Thus, $\delta(h)$ is a smooth functional on $P_{m_0}(M)$. Hence, it has an $\infty$-quasi-continuous modification (see [22]) which can be well defined outside of a slim set. Thus, we obtain a $\nu$-a.s. well-defined functional on the loop space $L_{m_0}(M)$. □

**4. The $L^p(\nu)$-convergence of divergence functionals.** For any $h\in H_0$ and $s<1$, the following functionals are well defined for $\mu$-a.s. $\gamma\in P_{m_0}(M)$:

(4.1) $\quad \delta_s(h)(\gamma)=\int_0^s(\dot h(\tau)+\tfrac{1}{2}\mathrm{Ric}_{r_x(\tau)}h(\tau), dx(\tau)), \qquad x=I^{-1}(\gamma).$

Since $\mu$ and $\nu$ are equivalent on $\mathcal{F}_s=\sigma(\gamma_x(\tau),\tau\leq s)$, $\delta_s(h)$ is also $\nu$-a.s. well defined on $L_{m_0}(M)$. The main technique part in Hsu's proof of the integration by parts formula on the loop space (see Proposition 4.1 in [15]) is to prove that, as $s \to 1$, the sequence $\{\delta_s(h)\}$ converges in $L^1(\nu)$ to a limit which belongs to $L^2(\nu)$. In view of this, for $\nu$-a.s. $\gamma \in L_{m_0}(M)$, Hsu defined $\delta(h)(\gamma)$ as the $L^1(\nu)$-limit of $\delta_s(h)(\gamma)$ and then proved that $\delta(h)$ is nothing else than the divergence functional appearing in an integration by parts formula on $L_{m_0}(M)$.

The purpose of this section is to prove that, for all $p \geq 1$, as $s \to 1$ the sequence $\{\delta_s(h)\}$ converges to $\widetilde{\delta(h)}$ in $L^p(\nu)$, where $\widetilde{\delta(h)}$ is an $\infty$-quasi-continuous modification of $\delta(h)$ constructed by the quasi-sure analysis principle which is $\nu$-a.s. well defined on $L_{m_0}(M)$ (see Corollary 3.5). Moreover, we prove that $\widetilde{\delta(h)}$ satisfies the Driver–Enchev–Stroock–Hsu integration by parts formula on the loop space.

THEOREM 4.1. *For any $p \geq 1$, there exists a constant $C_p > 0$ such that, for all $h \in H$, we have*

$$\|\widetilde{\delta(h)}\|_{L^p(\nu)} \leq C_p \|h\|_H. \tag{4.2}$$

*Moreover, for any $h \in H_0$,*

$$\|\delta_s(h) - \widetilde{\delta(h)}\|_{W^{r,p}(X)} \leq C \left( \int_s^1 |\dot{h}(t)|^2 \, dt \right)^{1/2}, \tag{4.3}$$

$$\|\delta_s(h) - \widetilde{\delta(h)}\|_{L^p(\nu)} \leq C \left( \int_s^1 |\dot{h}(t)|^2 \, dt \right)^{1/2}. \tag{4.4}$$

PROOF. By the Hölder inequality, we need only to prove Theorem 4.1 for $p = 2n$, $n \in \mathbb{N}$. By Theorem 2.1, there exist $r \in \mathbb{N}$ and $q > 1$ such that

$$\|\widetilde{\delta(h)}\|^p_{L^p(\nu)} \leq C \|(\widehat{\delta(h)})^p\|_{W^{r,q}(X)}.$$

It remains to prove that, for any $r \in \mathbb{N}$ and $q > 1$,

$$\|\nabla^r (\widehat{\delta(h)})^p\|_{L^q(\mu_0)} \leq C \|h\|_H^p. \tag{4.5}$$

By the Burkholder–Davis–Gundy inequality, we have

$$\|(\widehat{\delta(h)})^p\|_{L^q(\mu)} \leq C \|h\|_H^p. \tag{4.6}$$

Thus (4.5) holds for $r = 0$ and $q > 1$. Now we prove (4.5) for $r = 1$. By the Hölder inequality, we have

$$\begin{aligned}
\|\nabla(\widehat{\delta(h)})^p\|_{L^q(\mu_0)} &= \|p(\widehat{\delta(h)})^{p-1} \nabla \widehat{\delta(h)}\|_{L^q(\mu_0)} \\
&\leq p \|(\widehat{\delta(h)})^{p-1}\|_{L^{2q}(\mu_0)} \|\nabla \widehat{\delta(h)}\|_{L^{2q}(\mu_0)} \\
&\leq pC \|\widehat{\delta(h)}\|^{p-1}_{L^{2(p-1)q}(\mu_0)} \|h\|_H \qquad \text{[using (3.5) and (4.6)]} \\
&\leq pC \|h\|_H^p.
\end{aligned}$$

In general, for any $k \in \mathbb{N}$, we have
$$\nabla^{k+1} f^{2n} = 2n \nabla^k f^{2n-1} \otimes \nabla f + 2n f^{2n-1} \nabla^{k+1} f$$
$$= 2n(2n-1)\nabla^{k-1} f^{2n-2} \otimes \nabla f \otimes \nabla f$$
$$+ 2n(2n-1) f^{2n-2} \nabla^k f \otimes \nabla f + 2n f^{2n-1} \nabla^{k+1} f.$$

Thus, by induction and using the Hölder inequality together with (3.5) and (4.6), we can prove (4.5) for $p = 2n$ and all $r \in \mathbb{N}$.

For any $s < 1$, we can easily prove that $\delta_s(h) \circ I \in W^{\infty,\infty}(X)$. Hence $\delta_s(h) \in D^{\infty,\infty}(P_{m_0}(M))$. Let $\delta_s(h)^*$ be any $\infty$-quasi-continuous modification of $\delta_s(h)$. Then $\delta_s(h)^*$ is also well defined for $\nu$-a.s. $\gamma \in L_{m_0}(M)$ and $\delta_s(h) = \delta_s(h)^*$ holds $\nu$-a.s. on $L_{m_0}(M)$. By Theorem 2.1 and by the same argument used in the proof of (4.2) together with Proposition 3.4, we can prove that

$$\|\delta_{s_2}(h)^* - \delta_{s_1}(h)^*\|^p_{L^p(\nu)} = \|(\delta_{s_2}(h) - \delta_{s_1}(h))^*\|^p_{L^p(\nu)}$$
$$\leq C \| |\delta_{s_2}(h) - \delta_{s_1}(h)|^p \|_{W^{r,q}(X)}$$
(4.7)
$$= C \left\| \left| \int_{s_1}^{s_2} (\dot{h}(t) + J(x,t)h(t), dx(t)) \right|^p \right\|_{W^{r,q}}$$
$$\leq C \left[ \left( \int_{s_1}^{s_2} |\dot{h}(t)|^2 \, dt \right)^{p/2} + \left( \int_{s_1}^{s_2} |h(t)|^2 \, dt \right)^{p/2} \right],$$

where $C$ is a constant which only depends on $p$. Now for $h \in H_0$, we have $\int_s^1 |h(t)|^2 \, dt \leq \int_s^1 |\dot{h}(t)|^2 \, dt$. Taking $s_2 = 1$ in (4.7), we obtain (4.3) and (4.4). □

Now we introduce the gradient operator on $P_{m_0}(M)$ and $L_{m_0}(M)$. Let $F$ be a cylindrical functional on $P_{m_0}(M)$ [resp., $L_{m_0}(M)$] given by
$$F(\gamma) = f(\gamma_{s_1}, \ldots, \gamma_{s_k}),$$
where $f \in C^\infty(M^k)$, $0 = s_0 < s_1 < \cdots < s_k < s_{k+1} = 1$ is a finite partition of $[0,1]$, $k \in \mathbb{N}$. For $\mu$-a.s. $\gamma \in P_{m_0}(M)$ [resp., $\nu$-a.s. $\gamma \in L_{m_0}(M)$], we define the gradient $DF(\gamma)$ of $F$ as the unique element of $H$ (resp. $H_0$) such that, for any $h \in H$ (resp. $h \in H_0$),
$$\langle DF(\gamma), h \rangle_H = D_h F(\gamma).$$
Here
$$D_h F(\gamma) = \sum_{i=1}^k \langle \mathrm{grad}^{(i)} f(\gamma(s_i)), U_{s_i}(\gamma) h(s_i) \rangle, \qquad \nu\text{-a.s. } \gamma \in L_{m_0}(M),$$
and $\mathrm{grad}^{(i)} f(\gamma(s_i))$ denotes the gradient of $f$ with respect to the $i$th variable, $i = 1, \ldots, k$.

Below we prove that $\widetilde{\delta(h)}$ satisfies the Driver–Enchev–Stroock–Hsu integration by parts formula on the loop space.

THEOREM 4.2. *Let $F$, $G$ be two cylindrical functionals on $L_{m_0}(M)$, $h \in H_0$. Then*

$$\langle D_h F, G \rangle_{L^2(\nu)} = \langle F, D_h^* G \rangle_{L^2(\nu)},$$

*where*

$$D_h^* = -D_h + \widetilde{\delta(h)}.$$

PROOF. First, we suppose $h \in H_1 = \{h \in H_0 : \mathrm{supp}(h) \subset\subset (0,1)\}$. Then there exists some $\alpha \in [0,1)$ such that $h(\tau) = 0$, $\forall \tau \in [\alpha, 1]$. Without loss of generality, we may suppose that $F$ and $G$ are $\mathcal{F}_\alpha$-measurable. To simplify the notation, let $p_t(x) = p_t(m_0, x)$, $x \in M$. Note that

$$D_h p_{1-\alpha}(\gamma(\alpha)) = \mathrm{grad}\, p_{1-\alpha}(\gamma(\alpha)) U_\alpha(\gamma) h(\alpha) = 0.$$

By (1.1) and the integration by parts formula on the path space, we have

$$\int_{L_{m_0}(M)} D_h F(\gamma) G(\gamma)\, d\nu(\gamma)$$

$$= \int_{P_{m_0}(M)} D_h F(\gamma) G(\gamma) \frac{p_{1-\alpha}(\gamma(\alpha))}{p_1(m_0)}\, d\mu(\gamma)$$

$$= \int_{P_{m_0}(M)} F(\gamma) D_h^* \left( G(\gamma) \frac{p_{1-\alpha}(\gamma(\alpha))}{p_1(m_0)} \right) d\mu(\gamma)$$

$$= \int_{P_{m_0}(M)} F(\gamma) \left[ -D_h \left( G(\gamma) \frac{p_{1-\alpha}(\gamma(\alpha))}{p_1(m_0)} \right) \right.$$

$$\left. + \widetilde{\delta(h)} \left( G(\gamma) \frac{p_{1-\alpha}(\gamma(\alpha))}{p_1(m_0)} \right) \right] d\mu(\gamma)$$

$$= \int_{P_{m_0}(M)} F(\gamma) (-D_h G(\gamma) + \widetilde{\delta(h)} G(\gamma)) \frac{p_{1-\alpha}(\gamma(\alpha))}{p_1(m_0)}\, d\mu(\gamma)$$

$$- \int_{P_{m_0}(M)} F(\gamma) G(\gamma) D_h \left( \frac{p_{1-\alpha}(\gamma(\alpha))}{p_1(m_0)} \right) d\mu$$

$$= \int_{P_{m_0}(M)} F(\gamma) (-D_h G(\gamma) + \widetilde{\delta(h)} G(\gamma)) \frac{p_{1-\alpha}(\gamma(\alpha))}{p_1(m_0)}\, d\mu(\gamma)$$

$$= \int_{L_{m_0}(M)} F(\gamma) (-D_h G(\gamma) + \widetilde{\delta(h)} G(\gamma))\, d\nu(\gamma).$$

Hence Theorem 4.2 holds for $h \in H_1$.

Next, for any $h \in H_0$, since $H_1$ is dense in $H_0$, there exist $h_n \in H_1$ such that $\|h_n - h\|_{H_0} \to 0$. By (4.2) in Theorem 4.1, and using the fact that $\widetilde{\delta(h_n)} - \widetilde{\delta(h)} = \widetilde{\delta(h_n - h)}$, for any $p > 1$, we have $\|\widetilde{\delta(h_n)} - \widetilde{\delta(h)}\|_{L^p(\nu)} \leq C\|h_n - h\|_{H_0}$. In particular, $\lim_{n\to\infty} \|\widetilde{\delta(h_n)} - \widetilde{\delta(h)}\|_{L^2(\nu)} = 0$. Hence

$$\langle F, \widetilde{\delta(h_n)} G \rangle_{L^2(\nu)} \to \langle F, \widetilde{\delta(h)} G \rangle_{L^2(\nu)}, \qquad n \to \infty.$$

On the other hand, by the definition of $D_h$ and the Lebesgue dominated convergence theorem, we have

$$\langle D_{h_n} F, G \rangle_{L^2(\nu)} \to \langle D_h F, G \rangle_{L^2(\nu)}, \qquad n \to \infty,$$
$$\langle D_{h_n} G, F \rangle_{L^2(\nu)} \to \langle D_h G, F \rangle_{L^2(\nu)}, \qquad n \to \infty.$$

Note that $\langle D_{h_n} F, G \rangle_{L^2(\nu)} + \langle D_{h_n} G, F \rangle_{L^2(\nu)} = \langle F, \delta(h_n) G \rangle_{L^2(\nu)}$. Letting $n \to \infty$, we prove that Theorem 4.2 holds for all $h \in H_0$. □

REMARK 4.1. Combining Theorems 4.1 and 4.2 with Proposition 4.1 in [15] concerning the $L^1(\nu)$-convergence of $\delta_s(h)$ to $\delta(h)$, we conclude that $\widetilde{\delta(h)}$ is nothing else than the divergence function defined in [7] (in the case where $h \in \mathcal{C}^1 \cap H_0$) and [9, 15, 16]. In view of this and to simplify the notation, as we have done in the statement of Theorem 1.1, in the rest of this paper, we will use the notation $\delta(h)$ instead of $\widetilde{\delta(h)}$.

As a consequence of Theorem 4.2 and the $L^p(\nu)$-integrability of the divergence functional $\delta(h)$ (see Theorem 4.1), by a standard argument as used in [10] or [15, 16], we have the following result which allows us to introduce the first-order Sobolev spaces over the loop space.

THEOREM 4.3. *For all $p > 1$, $D$ is closable from $L^p(L_{m_0}(M), \nu)$ to $L^p(L_{m_0}(M), H, \nu)$.*

**5. Exponential integrability of divergence functionals.** In this section we prove a Fernique-type exponential integrability theorem for the divergence functional on the loop space.

THEOREM 5.1. *For all*

$$\lambda < \lambda_0 = \frac{1}{(2 + \|\operatorname{Ric}\|_\infty)\|h\|_H},$$

*we have*

(5.1) $$E_\nu[\exp(\lambda|\delta(h)|^2)] < +\infty,$$

*or, equivalently,*

$$\lim_{t\to\infty} \frac{1}{t^2} \log \nu(\{\gamma \in L_{m_0}(M) : |\delta(h)| > t\}) \leq -\frac{1}{(2 + \|\operatorname{Ric}\|_\infty)\|h\|_H},$$

*where $\|\operatorname{Ric}\|_\infty$ denotes the uniform bound of the Ricci curvature.*

To prove Theorem 5.1, we shall use the following lemma which provides us with a very useful tool to study the exponential integrability of some functionals with the pinned Wiener measure on loop spaces.

LEMMA 5.2.  *Let $f \in W^{n,n^2 p}(X)$, $n \in \mathbb{N}$, $p > 1$. Then*

$$(5.2) \qquad \|e^{f^2}\|_{W^{n,p}(X)} \leq C(n) \|e^{(1+\varepsilon)f^2}\|_p^{1/(1+\varepsilon)} \|f\|_{W^{n,n^2 p\varepsilon/(1+\varepsilon)}(X)}^n.$$

*Let $F \in D^{\infty,\infty}(P_{m_0}(M))$. Then there exist a constant $C > 0$ and a pair $(n,p) \in \mathbb{N} \times (1,+\infty)$ such that*

$$(5.3) \qquad E_\nu[e^{\widetilde{F}^2}] \leq C \|e^{(1+\varepsilon)F^2}\|_{L^p(\mu)}^{1/(1+\varepsilon)} \|F \circ I\|_{W^{n,n^2 p\varepsilon/(1+\varepsilon)}(X)}^n,$$

*where $\widetilde{F}$ denotes any $\infty$-quasi-continuous modification of $F$. Moreover, assuming that $M \subset \mathbb{R}^l$ is a Nash–Whitney embedding, then we can take $n = 2[\frac{l}{2}] + 2$ and $p > 1$ which can be arbitrarily close to $1$.*

PROOF.  By the chain rule, we have

$$\nabla^n e^{f^2} = e^{f^2} \sum_{r_1 + \cdots + r_n = n} c_{r_1,\ldots,r_n} \nabla^{r_1} f^2 \otimes \cdots \otimes \nabla^{r_n} f^2,$$

where $c_{r_1,\ldots,r_n}$ are some combinatorial constants which can be given explicitly, and the summation is taken over all $\{(r_1,\ldots,r_n) \in \mathbb{N}^n : 0 \leq r_1,\ldots,r_n \leq n, r_1 + \cdots + r_n = n\}$. Using the Hölder inequality, for any $\varepsilon > 0$, we have

$$\|\|\nabla^n e^{f^2}\|_{H^{\otimes n}}\|_p$$
$$\leq C(n) \|e^{f^2}\|_{p(1+\varepsilon)} \sum_{r_1 + \cdots + r_n = n} \|\nabla^{r_1} f^2\|_{np\varepsilon/(1+\varepsilon)} \cdots \|\nabla^{r_n} f^2\|_{np\varepsilon/(1+\varepsilon)}.$$

Moreover, for any $r \in \mathbb{N}$, there exist some constants $C_{j_1 \ldots j_r}$ such that

$$\nabla^r f^2 = \sum_{j_1 + \cdots + j_r = r} C_{j_1,\ldots,j_r} \nabla^{j_1} f \otimes \cdots \otimes \nabla^{j_r} f.$$

Hence

$$\|\nabla^r f^2\|_p \leq C(r) \sum_{j_1 + \cdots + j_r = r} \|\nabla^{j_1} f\|_{rp} \cdots \|\nabla^{j_r} f\|_{rp}$$
$$\leq C(r) \|f\|_{W^{r,rp}(X)}^r.$$

Thus,

$$\|\|\nabla^n e^{f^2}\|_{H^{\otimes n}}\|_p$$
$$\leq C(n) \|e^{(1+\varepsilon)f^2}\|_p^{1/(1+\varepsilon)}$$
$$\quad \times \sum_{r_1 + \cdots + r_n = n} \|f\|_{W^{r_1, r_1 np\varepsilon/(1+\varepsilon)}(X)}^{r_1} \cdots \|f\|_{W^{r_n, r_n np\varepsilon/(1+\varepsilon)}(X)}^{r_n}$$
$$\leq C(n) \|e^{(1+\varepsilon)f^2}\|_p^{1/(1+\varepsilon)} \|f\|_{W^{n,n^2 p\varepsilon/(1+\varepsilon)}(X)}^n.$$

Inequality (5.2) follows. Combining (5.2) with Theorem 2.2, we obtain (5.3). □

PROPOSITION 5.3. *For all*

$$\lambda < \lambda_0 = \frac{1}{(2 + \|\operatorname{Ric}\|_\infty)\|h\|_H},$$

*we have*

$$E_\mu[\exp(\lambda|\delta(h)|^2)] < +\infty.$$

PROOF. By random time changing, there exists a Brownian motion $\{B_t, t \in [0, \infty)\}$ which is adapted to the standard Brownian filtration $\mathcal{F}_s = \sigma(x(s), s \in [0, t])$ [here we allow $t \in [0, \infty)$] such that

$$\delta(h) \circ I(x) = B_T,$$

where

$$T = \left[\int_0^1 |\dot{h}(\tau) + \tfrac{1}{2}\operatorname{Ric}_{r_x(\tau)} h(\tau)|^2 \, d\tau\right]^{1/2} \leq \left(1 + \frac{\|\operatorname{Ric}\|_\infty}{2}\right)\|h\|_H.$$

By the refinement version of the well-known Fernique lemma ([11]; see also Theorem 3.3 in [19]), we have

$$E_\mu[\exp(\lambda\delta(h))] = E[\exp(\lambda\|B_T\|^2)]$$
$$\leq E\left[\exp\left(\lambda \sup_{s \in [0,(2+\|\operatorname{Ric}\|_\infty)\|h\|_H/2]} \|B_t\|^2\right)\right]$$
$$< +\infty,$$

provided that

$$\lambda(2 + \|\operatorname{Ric}\|_\infty)\frac{\|h\|_H}{2} < \frac{1}{2},$$

that is,

$$\lambda < \lambda_0 = \frac{1}{(2 + \|\operatorname{Ric}\|_\infty)\|h\|_H}. \quad \Box$$

PROOF OF THEOREM 5.1. Applying Lemma 5.2 to $F = \delta(h)$, for any $p > 1$ and $\varepsilon > 0$, there exists a constant $C > 0$ such that

$$E_\nu[\exp(\lambda|\delta(h)|^2)] \leq C\|\exp(\lambda(1+\varepsilon)|\widehat{\delta(h)}|^2)\|_{W^{n,p}(X)}^{1/(1+\varepsilon)} \|\widehat{\delta(h)}\|_{W^{n,n^2p\varepsilon/(1+\varepsilon)}(X)}^n,$$

where $n = 2[\frac{l}{2}] + 2$. By Proposition 5.3, we have

$$\|e^{\lambda(1+\varepsilon)p|\widehat{\delta(h)}|^2}\|_{W^{n,p}(X)} < +\infty,$$

provided that

$$\lambda < \frac{1}{p(1+\varepsilon)(2+\|\operatorname{Ric}\|_\infty)\|h\|_H}.$$

On the other hand, Theorem 3.3 says

$$\|\widehat{\delta(h)}\|_{W^{n,n^2 p\varepsilon/(1+\varepsilon)}(X)} < C\|h\|_H.$$

Thus, for all $\lambda < [p(1+\varepsilon)(2+\|\operatorname{Ric}\|_\infty)\|h\|_H]^{-1}$, we have $E_\nu[\exp(\lambda|\delta(h)|^2)] < +\infty$. Since we can choose $p$ arbitrarily close to 1 and $\varepsilon$ arbitrary close to 0, we deduce the desired inequality (5.1) for all $\lambda < \lambda_0 = [(2+\|\operatorname{Ric}\|_\infty)\|h\|_H]^{-1}$. □

As a consequence of Theorem 5.1, we have the following result.

THEOREM 5.4. *For any $\lambda > 0$ and any $h \in H_0$, we have*

$$E_\nu[\exp(\lambda|\delta(h)|)] < +\infty.$$

PROOF. By Theorem 5.1, we have

$$\lim_{t\to\infty} \frac{1}{t^2} \log \nu(\{\gamma \in L_{m_0}(M) : |\delta(h)| > t\}) \leq -\frac{1}{(2+\|\operatorname{Ric}\|_\infty)\|h\|_H}.$$

This yields that

$$\lim_{t\to\infty} \frac{1}{t} \log \nu(\{\gamma \in L_{m_0}(M) : |\delta(h)| > t\}) = -\infty.$$

Theorem 5.4 follows. □

**6. Smoothness of Driver's flow on path spaces.** Recall that by [6, 8, 14], for all $h \in H$, the vector field $D_h$ generates a global flow $\{\Phi_t, t \in \mathbb{R}\}$, the so-called Driver flow on $P_{m_0}(M)$, such that, for $\mu$-a.s. $\gamma \in P_{m_0}(M)$,

(6.1)
$$\dot{\Phi}_t(\gamma) = D_h(\Phi_t(\gamma)),$$
$$\Phi_0(\gamma) = \gamma.$$

In the case where $h \in \mathcal{C}^1 \cap H_0$ is a Lipschitz Cameron–Martin vector, Driver [7] constructed the flow generated by $D_h$ on the loop space $L_{m_0}(M)$ by using the technique of enlargement of filtration. In [9], Enchev and Stroock gave another approach to construct the flow corresponding to $D_h$ for all $h \in H_0$. As we have pointed out before, their approaches relied on the gradient estimate and the Hessian estimate of the logarithm of the heat kernel on compact Riemannian manifold. See also Section 9.

In this section, without using any heat kernel estimate, we construct the Driver flow of $D_h$ on the loop space for all $h \in H_0$ through an $\infty$-quasi-continuous modification of the corresponding flow on the path space. To this end, we shall first prove that, for any $h \in H_0$, the Driver flow generated by $D_h$ is a smooth transform on the path space $P_{m_0}^{2m,\alpha}(M)$ in the sense of the Malliavin calculus, where for any $m \in \mathbb{N}$, $m \geq 2$, $\alpha \in (\frac{1}{2m}, \frac{1}{2})$,

$$P_{m_0}^{2m,\alpha}(M) = \left\{ \gamma \in P_{m_0}(M) : \int_0^1 \int_0^1 \frac{d(\gamma(t), \gamma(s))^{2m}}{|t-s|^{1+2m\alpha}} \, dt \, ds < \infty \right\}.$$

By [22] and the references therein, the Wiener measure $\nu$ is supported on $P_{m_0}^{2m,\alpha}(M)$ and the Sobolev spaces theory on $P_{m_0}^{2m,\alpha}(M)$ is the same as (i.e., quasi-homeomorphic to) the one on $P_{m_0}(M)$. Moreover, $P_{m_0}^{2m,\alpha}(M)$ is an $M$-type 2 Banach manifold modeled on $X^{2m,\alpha}$ whose norm $\|\cdot\|_{2m,\alpha}$ is smooth in $X^{2m,\alpha} \setminus \{0\}$, where

$$X^{2m,\alpha} := \left\{ x \in X : \int_0^1 \int_0^1 \frac{\|x(t) - x(s)\|_{\mathbb{R}^d}^{2m}}{|t-s|^{1+2m\alpha}} \, dt \, ds < \infty \right\},$$

on which we consider the fractional Hölder norm $\|\cdot\|_{2m,\alpha}$ given by

$$\|x\|_{2m,\alpha} = \left( \int_0^1 \int_0^1 \frac{\|x(t) - x(s)\|_{\mathbb{R}^d}^{2m}}{|t-s|^{1+2m\alpha}} \, dt \, ds \right)^{1/(2m)}.$$

THEOREM 6.1. *For any $h \in H$, Driver's flow $\Phi_t : P_{m_0}^{2m,\alpha}(M) \to P_{m_0}^{2m,\alpha}(M)$ is a smooth mapping in the sense of the Malliavin calculus, that is,*

$$\Phi_t \in D^{\infty,\infty}(P_{m_0}^{2m,\alpha}(M), P_{m_0}^{2m,\alpha}(M)).$$

To prove this theorem, let us first introduce the set $\mathcal{SM}(\mathbb{R}^d)$ of all $\mathbb{R}^d$-valued semimartingales with the Doob–Meyer decomposition $\xi(s) = \int_0^s O(r) \, dx(r) + \int_0^s A(r) \, dr$, where $\{(O(s), A(s)), s \in [0,1]\}$ is an adapted $M(d,d) \times \mathbb{R}^d$-valued process such that $\|\|\xi\|\|_2^2 := E[\sup_{s \in [0,1]} \|O(s)\|_{M(d,d)}^2] + E[\int_0^1 \|A(s)\|_{\mathbb{R}^d}^2 \, ds]$ is finite. By definition, we have

$$\left\| \|\xi(\cdot)\|_{2m,\alpha} \right\|_p = \left\| \left\| \int_0^\cdot O(r) \, dx(r) + \int_0^\cdot A(r) \, dr \right\|_{2m,\alpha} \right\|_p$$

$$\leq \left\| \left\| \int_0^\cdot O(r) \, dx(r) \right\|_{2m,\alpha} \right\|_p + \left\| \left\| \int_0^\cdot A(r) \, dr \right\|_{2m,\alpha} \right\|_p.$$

Using the Burkholder–Davis–Gundy inequality, we have

$$E\left[ \left\| \int_0^\cdot O(r) \, dx(r) \right\|_{2m,\alpha}^p \right] \leq C_p E\left[ \left\| \int_0^\cdot |O(r)|^2 \, dr \right\|_{2m,\alpha}^{p/2} \right]$$

$$\leq C_p E\left[\sup_{r\in[0,1]} \|O(r)\|^{2\cdot p/2} \left\|\int_0^\cdot dr\right\|_{2m,\alpha}^{p/2}\right]$$

$$\leq C(p,m,\alpha) E\left[\sup_{r\in[0,1]} \|O(r)\|^p\right].$$

Using the Cauchy–Schwarz inequality, we have

$$\left|\int_{s_1}^{s_2} A(r)\,dr\right|^{2m} \leq \left[\int_0^1 |A(r)|^2\,dr\right]^m |s_1-s_2|^m.$$

Thus

$$\int_0^1 \int_0^1 \frac{|\int_{s_1}^{s_2} A(r)\,dr|^{2m}}{|s_1-s_2|^{1+2m\alpha}}\,ds_1\,ds_2 \leq \left[\int_0^1 |A(r)|^2\,dr\right]^m \int_0^1\int_0^1 \frac{ds_1\,ds_2}{|s_1-s_2|^{1+2m\alpha-m}}$$

$$\leq C(m,\alpha)\left[\int_0^1 |A(r)|^2\,dr\right]^m.$$

This yields that

$$\left\|\left\|\int_0^\cdot A(r)\,dr\right\|_{2m,\alpha}\right\|_p \leq [C(m,\alpha)]^{1/2m}\left\{E\left[\int_0^1 |A(r)|^2\,dr\right]^{pm\times 1/2m}\right\}^{1/p}$$

$$\leq C(m,\alpha)\left\{E\left[\int_0^1 |A(r)|^2\,dr\right]^{p/2}\right\}^{1/p}.$$

Therefore,

$$\|\|\xi\|_{2m,\alpha}\|_p \leq C(m,p,\alpha)\left\{E\left[\sup_{r\in[0,1]} |O(r)|^p\right]\right\}^{1/p}$$

$$+ C(m,p,\alpha)\left\{E\left[\int_0^1 |A(r)|^2\,dr\right]^{p/2}\right\}^{1/p}.$$

Hence for $m \geq 2$, $\alpha \in (\frac{1}{2m}, \frac{1}{2})$ and $p \geq 1$, we have

$$E[\|\xi\|_{2m,\alpha}^p] \leq C(m,p,\alpha) E\left[\sup_{r\in[0,1]} \|O(r)\|^p\right] + C(m,p,\alpha) E\left[\int_0^1 |A(r)|^2\,dr\right]^{p/2}.$$

For any $p \geq 1$, define the norm $\|\|\cdot\|\|_p$ on $\mathcal{SM}(\mathbb{R}^d)$ as

$$\|\|\xi\|\|_p^p := E\left[\sup_{r\in[0,1]} \|O(r)\|^p\right] + E\left[\int_0^1 |A(r)|^2\,dr\right]^{p/2}.$$

Then

(6.2) $$E[\|\xi\|_{2m,\alpha}^p] \leq C(m,p,\alpha) \|\|\xi\|\|_p^p.$$

PROOF OF THEOREM 6.1. Set $\xi_t = I^{-1} \circ \Phi_t \circ I$. By [22], we have $I \in D^{\infty,\infty}(X^{2m,\alpha}, P_{m_0}^{2m,\alpha}(M))$ and $I^{-1} \in D^{\infty,\infty}(P_{m_0}^{2m,\alpha}(M), X^{2m,\alpha})$. Hence it remains to prove $\xi_t \in W^{\infty,\infty}(X^{2m,\alpha}, X^{2m,\alpha})$. By [6, 14], $\xi_t$ satisfies the following ODE on $X^{2m,\alpha}$, where the stochastic integral is taken in the sense of Itô:

$$\text{(6.3)} \quad \frac{\partial}{\partial t}\xi_t(s) = h(s) + \frac{1}{2}\int_0^s \text{Ric}_{U_t(\tau)}(h(\tau))\,d\tau + \int_0^s q_h(t,\tau)\,d\xi_t(\tau),$$
$$\xi_0(x) = x.$$

Here $U_t(s)$ is the stochastic parallel transport along $\gamma_t(s) = \Phi_t(\gamma)(s)$ and $q_h(t,s)$ is given by the following Stratonovich stochastic integral:

$$q_h(t,s) = \int_0^s \Omega_{U_t(\tau)}(h(\tau) \circ d\xi_t(\tau)).$$

Using the Picard iteration and by a similar argument as used in [6] and [15], we can prove that, for any $k \in H$, $D_k\xi_t(s)$ exists for all $s \in [0,1]$ and all $t \in \mathbb{R}$. Moreover, $D_k\xi_t$ satisfies the following equation:

$$\frac{\partial}{\partial t}D_k\xi_t(s) = \int_0^s q_h(t,\tau)\,dD_k\xi_t(\tau)$$
$$\text{(6.4)} \qquad\qquad + \frac{1}{2}\int_0^s D_k\text{Ric}_{U_t(\tau)}(h(\tau))\,d\tau + \int_0^s D_k q_h(t,\tau)\,d\xi_t(\tau),$$
$$D_k\xi_0(s) = k(s).$$

For any $T > 0$, $p > 1$, we can prove

$$\text{(6.5)} \qquad \sup_{t \in [-T,T]} E[\|\|D\xi_t(\cdot)\|_H\|_{2m,\alpha}^p] < \infty.$$

Indeed, let $\{\xi_t^n(s), s \in [0,1], t \in [-T,T]\}$ be given as in [6] and [14]. Let $(O_t^n(s), A_t^n(s))$ be the Doob–Meyer decomposition of $\xi_t^n$. Let $\gamma_t^n = I(\xi_t^n)$, $U_t^n = U(\gamma_t^n)$ and $q_h^n(t,s) = \int_0^s \Omega_{U_t^n(r)}(h(r) \circ d\xi_t^n(r))$, $s \in [0,1]$, $t \in [-T,T]$. Then for any fixed $s \in [0,1]$ and $t \in [-T,T]$, it is easy to see that $(O_t^n(s), A_t^n(s)) \in W^{\infty,\infty}(X, O(d) \times \mathbb{R}^d)$. Similarly to [14], we can easily prove that

$$\text{(6.6)} \qquad \sup_{t \in [-T,T]} E[\|\|\xi_t^n - \xi_t^{n-1}\|\|^p] \leq \frac{(CT)^n}{n!},$$

where

$$E[\|\|\xi_t^n - \xi_t^{n-1}\|\|^p] := E\left[\sup_{s \in [0,1]} \|O_t^n(s) - O_t^{n-1}(s)\|_{M(d,d)}^p\right]$$
$$+ E\left[\int_0^1 \|A_t^n(s) - A_t^{n-1}(s)\|^2\,ds\right]^{p/2}.$$

Moreover, for all $k \in H$, we can easily show that

$$D_k O_t^n(s) = O_t^n(s) \int_0^t [O_u^n(s)]^{-1} D_k q_h^{n-1}(u,s) O_u^n(s)\, du,$$

$$D_k A_t^n(s) = O_t^n(s) \int_0^t [O_u^n(s)]^{-1} D_k q_h^{n-1}(u,s) A_u^n(s)\, du$$
$$+ \tfrac{1}{2} O_t^n(s) \int_0^t [O_u^n(s)]^{-1} D_k \operatorname{Ric}_{U_u^{n-1}(s)}(h(s))\, du.$$

Putting $k \in H$ such that $\dot{k} = \mathbf{1}_{[\tau,1]} e_\alpha$ into the above formulas, we obtain the explicit expressions of the Malliavin derivatives $D_\tau^\alpha O_t^n(s)$ and $D_\tau^\alpha A_t^n(s)$. Set

$$E[\|\!|D_\tau^\alpha \xi_t^n|\!\|^p] := E[\|O_t^n(\tau)\|_{M(d,d)}^p] + E\left[\sup_{s \in [0,1]} \|D_\tau^\alpha O_t^n(s)\|_{M(d,d)}^p\right]$$
$$+ E\left[\int_0^1 \|D_\tau^\alpha A_t^n(s)\|_{\mathbb{R}^d}^2\, ds\right]^{p/2}.$$

By standard argument and the Burkholder–Davis–Gundy inequality, and using the fact that $O_t^n(s) \in O(d)$ and $\sup_{n \in \mathbb{N}} \|A_t^n(s)\| \leq C(1 + |\dot{h}(s)|)$ (see [14]), it is straightforward to prove that

$$E\left[\sup_{s \in [0,1]} \|D_\tau^\alpha U_t^n(s)\|^p\right] \leq cE[\|\!|D_\tau^\alpha \xi_t^n|\!\|^p],$$

$$E\left[\sup_{s \in [0,1]} \|D_\tau^\alpha q_h^n(t,s)\|^p\right] \leq cE[\|\!|D_\tau^\alpha \xi_t^n|\!\|^p],$$

$$E\left[\sup_{s \in [0,1]} \|D_\tau^\alpha U_t^n(s) - D_\tau^\alpha U_t^{n-1}(s)\|^p\right]$$
$$\leq cE[\|\!|D_\tau^\alpha \xi_t^n - D_\tau^\alpha \xi_t^{n-1}|\!\|^p]$$
$$+ c\left\{E\left[\sup_{s \in [0,1]} \|D_\tau^\alpha U_t^n(s)\|^{2p}\right]\right\}^{1/2} \{E[\|\!|\xi_t^n - \xi_t^{n-1}|\!\|^{2p}]\}^{1/2},$$

$$E\left[\sup_{s \in [0,1]} \|D_\tau^\alpha q_h^n(t,s) - D_\tau^\alpha q_h^{n-1}(t,s)\|^p\right]$$
$$\leq cE[\|\!|D_\tau^\alpha \xi_t^n - D_\tau^\alpha \xi_t^{n-1}|\!\|^p]$$
$$+ c\{E[\|\!|D_\tau^\alpha \xi_t^{n-1}|\!\|^{2p}]\}^{1/2} \{E[\|\!|\xi_t^{n-1} - \xi_t^{n-2}|\!\|^{2p}]\}^{1/2}.$$

From the above inequalities, we can deduce that

(6.7) $$\sup_{n \in \mathbb{N}} \sup_{t \in [-T,T]} \sup_{\tau \in [0,1]} E[\|\!|D_\tau^\alpha \xi_t^n|\!\|^p] < c_1 e^{c_2 T},$$

$$E\left[\sup_{s\in[0,1]} \|D^\alpha_\tau O^n_t - D^\alpha_\tau O^{n-1}_t(s)\|^p\right]$$

(6.8)
$$\leq R^1_{\tau,\alpha}(T,n) + c\int_0^t E[\||D^\alpha_\tau \xi^n_u - D^\alpha_\tau \xi^{n-1}_u\||^p]\, du,$$

$$E\left[\left|\int_0^1 \|D^\alpha_\tau A^n_t - D^\alpha_\tau A^{n-1}_t(s)\|^2\, ds\right|^{p/2}\right]$$

(6.9)
$$\leq R^2_{\tau,\alpha}(T,n) + c\int_0^t E[\||D^\alpha_\tau \xi^n_u - D^\alpha_\tau \xi^{n-1}_u\||^p]\, du,$$

where for $i=1,2$, there exists a constant $C(T,p)$ such that

(6.10)
$$R^i_{\tau,\alpha}(T,n) \leq C(T,p)\left\{\sup_{t\in[-T,T]}\sup_{\tau\in[0,1]} E[\||D^\alpha_\tau \xi^n_t\||^{2p}]\right\}^{1/2}$$
$$\times \left\{\sup_{t\in[-T,T]} E[\||\xi^n_t - \xi^{n-1}_t\||^{2p}]\right\}^{1/2}.$$

From (6.6), (6.7) and (6.10), $R_{\tau,\alpha}(T,n) = R^1_{\tau,\alpha}(T,n) + R^2_{\tau,\alpha}(T,n)$ tends to zero as $n$ tends to infinity. hand, from (6.8), (6.9) and the Gronwall inequality, we have

(6.11)
$$\sup_{t\in[-T,T]}\sup_{\tau\in[0,1]} E[\||D^\alpha_\tau \xi^n_t - D^\alpha_\tau \xi^{n-1}_t\||^p] \leq R_{\tau,\alpha}(T,n)e^{cT}.$$

This implies that $\{\xi^n_t, t\in[-T,T]\}$ converges uniformly in $W^{1,\infty}(X,(\mathcal{S}M(\mathbb{R}^d), \|\|\cdot\|\|_p))$ and hence by (6.2) it converges uniformly in $W^{1,\infty}(X,(X^{2m,\alpha}, \|\cdot\|_{2m,\alpha}))$. Moreover, we deduce (6.4) [resp., the inequality (6.5)] from the corresponding equation for $D_k\xi^n_t$ [resp., the inequality (6.7)]. In general, by induction and repeating the same argument as above, we can prove that, for all $k_1,\ldots,k_r \in H$, $D_{k_1,\ldots,k_r}\xi_t(s)$ exists for all $s\in[0,1]$ and all $t\in\mathbb{R}$. Moreover, for any $T>0$, $p>1$, using the Burkholder–Davis–Gundy inequality, we have

(6.12)
$$\sup_{t\in[-T,T]} E[\||\|D^r\xi_t(\cdot)\|_{H^{\otimes r}}\|^p_{2m,\alpha}] < \infty.$$

This completes the proof of Theorem 6.1 concerning the smoothness of the Driver flow on the path space $P^{2m,\alpha}_{m_0}(M)$. To save the length of the paper, we omit the details of the proofs of the four inequalities listed before (6.7). The reader who is interested in the details of the proof is referred to [3] (for the case where $h\in\mathcal{C}^1\cap H$ is a Lipschitz Cameron–Martin vector) and [18] as well as [21] where the author proved that the Driver flow $\Phi_t$ is a smooth transform on the path space $P_{m_0}(M)$. □

We will make use of the following Kolmogorov criterion for $\infty$-quasi-continuous modification of a family of $M$-type 2 Banach space $E$-valued functionals. When $E = \mathbb{R}$, it is due to [28]. See also [26].

THEOREM 6.2. *Let $\{X(t), t \in [-T,T]\}$ be a family of $M$-type 2 Banach space $E$-valued functional. Suppose that, for all $p \geq 2$, $r \in \mathbb{N}$, there exist constants $c, \varepsilon > 0$ and an even number $\beta$ such that:*

(i) $X(t) \in W^{r,p}(X, E)$;
(ii) *for all* $(s,t) \in [-T,T] \times [-T,T]$ $\|X(t) - X(s)\|_E^\beta \in W^{r,p}(X)$;
(iii) *for all* $(s,t) \in [-T,T] \times [-T,T]$, *we have*

$$\|\|X(t) - X(s)\|_E^\beta\|_{W^{r,p}(X)} \leq c|t-s|^{1+\varepsilon}.$$

*Then there exists a version of the process $\{X(t), t \in [-T,T]\}$ which is $\infty$-quasi-continuous for each $t \in [-T,T]$ and which has continuous paths.*

PROOF. Since $E$ is an $M$-type 2 Banach space, the norm $\phi(x) = \|x\|_E$ is smooth in $E \setminus \{0\}$ and, for all $k \in \mathbb{N}$, there exists $M_k$ such that $\sup_{\|x\|_E = 1} \|\nabla^k \times \phi\|(x) \leq M_k < \infty$. Thus, the Chebyshev-type inequality of $(r,p)$-capacity for $E$-valued functionals holds; see, for example, [26]. Hence, for any given $r \in \mathbb{N}$, $p > 1$, $\varepsilon > 0$, we have

$$C_{r,p}(\{x \in X : \|X(t) - X(s)\|_E > \varepsilon\}) \leq \frac{1}{\varepsilon^\beta} \|\|X(t) - X(s)\|_E^\beta\|_{r,p}.$$

Therefore, Theorem 6.2 can be proved by the same argument as used in the proof of Theorem 3.1 in [28]. □

THEOREM 6.3. *For all $p \geq 2$ and $r \in \mathbb{N}$, there exist constants $c, \varepsilon > 0$ such that:*

(i) $\xi_t \in W^{r,p}(X^{2m,\alpha}, X^{2m,\alpha})$;
(ii) $\|\xi_t - \xi_s\|_{2m,\alpha}^\beta \in W^{r,p}(X^{2m,\alpha})$ *for all* $(s,t) \in [-T,T] \times [-T,T]$;
(iii) *for all* $(s,t) \in [-T,T] \times [-T,T]$, *we have*

(6.13) $$\|\|\xi_t(\cdot) - \xi_s(\cdot)\|_{2m,\alpha}^\beta\|_{r,p} \leq c|t-s|^{1+\varepsilon}.$$

PROOF. By Theorem 6.1, for all $t \in [-T,T]$, we have $\xi_t \in W^{r,p}(X^{2m,\alpha}, X^{2m,\alpha})$. By the chain rule, for all even number $\beta > 0$ and $(s,t) \in [-T,T] \times [-T,T]$, we can prove $\|\xi_t - \xi_s\|_{2m,\alpha}^\beta \in W^{r,p}(X^{2m,\alpha})$. By Lemma 4.1 in [28], for all $(s,t) \in [-T,T] \times [-T,T]$, we have

(6.14) $$\|\|\xi_t(\cdot) - \xi_s(\cdot)\|_{2m,\alpha}^n\|_{2r,p} \leq C(n,p,r) \|\|\xi_t(\cdot) - \xi_s(\cdot)\|_{2m,\alpha}\|_{2r,4r^2p}^{2r}$$
$$\times \max_{0 \leq k \leq n} [E\|\xi_t(\cdot) - \xi_s(\cdot)\|_{2m,\alpha}^{(n-k)2rp}]^{1/2rp}.$$

On the other hand, using (6.3), (6.4) and inequalities (6.5) and (6.12), by the Hölder inequality and the Burkholder–Davis–Gundy inequality, we can verify that, for all $p \geq 1$ and $r \in \mathbb{N}$,

$$\|\|\xi_t(\cdot) - \xi_s(\cdot)\|_{2m,\alpha}\|_{2p}^{2p} \leq C(p,T)|t-s|^p, \tag{6.15}$$

$$\|\|\xi_t(\cdot) - \xi_s(\cdot)\|_{2m,\alpha}\|_{r,2p}^{2p} \leq C(p,T)|t-s|^p. \tag{6.16}$$

From (6.14), (6.15) and (6.16), we deduce that, for $n(p,r) = n$ large enough and for some constant $\varepsilon > 0$, we have

$$\|\|\xi_t(\cdot) - \xi_s(\cdot)\|_{2m,\alpha}^{n(p,r)}\|_{2r,p} \leq c|t-s|^{1+\varepsilon}.$$

The proof of Theorem 6.3 is complete. $\square$

Combining Theorem 6.3 with the Kolmogorov criterion (Theorem 6.2), we have the following:

THEOREM 6.4. *For all $T > 0$, there exists a version of the Driver flow $\{\Phi_t, t \in [-T,T]\}$ on the path space $P_{m_0}^{2m,\alpha}(M)$ which is $\infty$-quasi-continuous for each $t \in [-T,T]$ and which has continuous trajectory on $t \in [-T,T]$.*

PROOF OF THEOREM 1.4. By Theorem 6.3, for any $T > 0$, there exists a slim subset $S_T$ of the path space $P_{m_0}^{2m,\alpha}(M)$ such that an $\infty$-quasi-continuous modification (denoted by $\{\widetilde{\Phi}_t, t \in [-T,T]\}$) of $\{\Phi_t, t \in [-T,T]\}$ can be well defined for all $\gamma \in P_{m_0}^{2m,\alpha}(M) \setminus S_T$. Taking $T_n = 2^n$, we deduce that there exists a common slim set $S_\infty = \bigcup_{n \in \mathbb{N}} S_{T_n}$ such that $\{\widetilde{\Phi}_t, t \in \mathbb{R}\}$ can be well defined for all $\gamma \in P_{m_0}^{2m,\alpha}(M) \setminus S_\infty$. Thus, $\{\widetilde{\Phi}_t, t \in \mathbb{R}\}$ can be $\nu$-a.s. well defined on the loop space $L_{m_0}^{2m,\alpha}(M) = L_{m_0}(M) \cap P_{m_0}^{2m,\alpha}(M)$ and hence is $\nu$-a.s. well defined on $L_{m_0}(M)$. $\square$

In the rest of this paper, we fix such an $\infty$-quasi-continuous version $\{\widetilde{\Phi}_t, t \in \mathbb{R}\}$ of $\{\Phi_t, t \in \mathbb{R}\}$. To end this section, let us mention the following remark.

REMARK 6.1. Note that the Driver flow $\{\Phi_t, t \in \mathbb{R}\}$ is $\mathcal{F}_s/\mathcal{F}_s$-measurable for all $s \in [0,1]$; see [6, 8, 9]. Thus, the Kolmogorov criterion yields that $\{\widetilde{\Phi}_t, t \in \mathbb{R}\}$ is again $\mathcal{F}_s/\mathcal{F}_s$-measurable for all $s \in [0,1]$. Moreover, using the same argument as used in the proof of the Sobolev norm and the capacity comparison theorems between the Wiener space and the path space via the Itô map (see [22]), we can prove that, for any $F \in W^{r,p}(X)$ and any subset $A \subset X$, ant $r \in \mathbb{N}$ and $p > 1$,

$$\alpha_1 \|F \circ \xi_t\|_{r/2,p-\varepsilon} \leq \|F\|_{r,p} \leq \alpha_2 \|F \circ \xi_t\|_{2r,p+\varepsilon}, \tag{6.17}$$

$$\alpha_1 C_{r/2,p-\varepsilon}(\widetilde{\xi}_t^{-1}(A)) \leq C_{r,p}(A) \leq \alpha_2 C_{2r,p+\varepsilon}(\widetilde{\xi}_t^{-1}(A)), \tag{6.18}$$

where $\|\cdot\|_{r,p}$ (resp., $C_{r,p}$) denotes the $(r,p)$-Sobolev norm [resp., $(r,p)$-capacity] on the Wiener space $X$, $\{\widetilde{\xi}_t, t \in \mathbb{R}\}$ denotes any $\infty$-quasi-continuous modification of the pull-back of the Driver flow $\xi_t = I^{-1} \circ \Phi_t \circ I$ and $\alpha_1$ and $\alpha_2$ are two constants which depend only on $r$, $p$ and the uniform bounds of the Riemannian curvature and the Ricci curvature as well as their higher-order covariant derivatives. This yields that the flow property $\widetilde{\xi}_s \circ \widetilde{\xi}_t = \widetilde{\xi}_{t+s}$ holds quasi-surely on $X$. Since $\widetilde{\Phi}_t = \widetilde{I} \circ \widetilde{\xi}_t \circ \widetilde{I}^{-1}$, we get the flow property $\widetilde{\Phi}_s \circ \widetilde{\Phi}_t = \widetilde{\Phi}_{t+s}$ quasi-surely on $P_{m_0}(M)$. As explained in [22], the proof of the above inequalities (6.17) and (6.18) are based on the Meyer inequality on the Wiener space. However, we still do not know whether or not the Meyer inequality holds on the path space over a compact Riemannian manifold. Thus, we do not know whether or not the corresponding Sobolev norms (resp., capacities) comparison inequalities hold on path spaces if one replaces $\xi_t$ by $\Phi_t$.

**7. Cameron–Martin theorem on loop spaces.** In this section we will first construct the Driver flow on the loop space through the corresponding flow on the path space. Combining this and Theorem 5.4 together with the Cruzeiro lemma, we will give an alternative proof to the Cameron–Martin theorem on loop spaces established earlier by Driver [7] and Enchev and Stroock [9] by Doob's $h$-processes approach and the short time upper bound estimates of the gradient and the Hessian of logarithm of the heat kernels.

Our first result in this section is the following theorem:

THEOREM 7.1. *Let $h \in H_0$. Then $\widetilde{\Phi}_t(L_{m_0}(M)) \subset L_{m_0}(M)$. Moreover, for $\nu$-a.s. $\gamma \in L_{m_0}(M)$, we have*

$$\dot{\widetilde{\Phi}}_t(\gamma) = D_h(\widetilde{\Phi}_t(\gamma)),$$

(7.1)

$$\widetilde{\Phi}_t(\gamma) = \gamma.$$

In view of Theorem 7.1, we regard $\widetilde{\Phi}_t$ as the flow on $L_{m_0}(M)$ generated by $D_h$.

PROOF OF THEOREM 7.1. Since $\widetilde{\Phi}_t = \Phi_t$ holds quasi-surely on $P_{m_0}(M)$, the flow equation (6.1) is verified. It remains to show that $\widetilde{\Phi}_t(L_{m_0}(M)) \subset L_{m_0}(M)$, that is,

(7.2) $$\widetilde{\Phi}_t(\gamma)(1) = m_0.$$

To this end, we use the same argument as in [7]. Indeed, for $\nu$-a.s. $\gamma \in L_{m_0}(M)$ and any $s < 1$, by the flow property $\widetilde{\Phi}_t$ on $L_{m_0}(M)$ we have

$$d(\widetilde{\Phi}_t(\gamma)(s), \widetilde{\Phi}_0(\gamma)(s)) \leq \int_0^t \left|\frac{d}{d\tau}\widetilde{\Phi}_\tau(\gamma)(s)\right| d\tau$$

$$= \int_0^t |U(\widetilde{\Phi}_\tau(\gamma))(s)h(s)|\,d\tau.$$

Now $U(\widetilde{\Phi}_\tau(\gamma))(s)$ is an isometry from $\mathbb{R}^d$ to $T_{\widetilde{\Phi}_\tau(\gamma)(s)}M$. Thus for $\nu$-a.s. $\gamma \in L_{m_0}(M)$,

(7.3) $$d(\widetilde{\Phi}_t(\gamma)(s), \widetilde{\Phi}_0(\gamma)(s)) \leq |h(s)|.$$

Let

$$\Sigma = \{\gamma \in L_{m_0}(M) : d(\widetilde{\Phi}_t(\gamma)(s), \widetilde{\Phi}_0(\gamma)(s)) \leq |h(s)| \ \forall s < 1\}.$$

Since both sides of the inequality in (7.3) are continuous, we have

(7.4) $$\nu(\Sigma) = 1.$$

Taking $\gamma \in \Sigma$ and letting $s \to 1$, we have

$$\lim_{s \to 1} d(\widetilde{\Phi}_t(\gamma)(s), \gamma(s)) \leq \lim_{s \to 1} |h(s)| = 0.$$

By the continuity of $s \to \widetilde{\Phi}_t(s)$, we prove (7.2) for $\gamma \in \Sigma$, which differs from $L_{m_0}(M)$ up to a $\nu$-negligible subset, compare (7.4). $\square$

The following lemma is due to Cruzeiro [4] and is a very useful tool to study the quasi-invariance of a probability measure under the action of certain flows.

LEMMA 7.2 ([4]). *Let $(\Omega, \mathcal{F}, \{\mathcal{F}_s\}, P)$ be a complete filtered probability space and let $\Phi = \{\phi_t\}$ be a flow (i.e., a one-parameter group of measurable transformations) on $\Omega$. Suppose that there exists the divergence $\mathrm{div}(\Phi) \in L^1(\Omega, P)$ such that, for all $f \in \mathcal{C} \subset L^\infty(\Omega, P)$, where $\mathcal{C}$ is a dense subset of $L^\infty(\Omega, P)$, we have*

(7.5) $$\left\{\frac{d}{dt}E[f(\phi_t)]\right\}\bigg|_{t=0} = E[f\,\mathrm{div}(\Phi)].$$

*Moreover, assume that there exists a $\lambda > 0$ such that*

(7.6) $$E[e^{\lambda|\mathrm{div}(\Phi)|}] < +\infty.$$

*Then $(\phi_t)_* P$ is absolutely continuous with respect to $P$. Denote*

$$K_t = \frac{d(\phi_t)_* P}{dP}.$$

*Then*

$$K_t = \exp\left(\int_0^t \mathrm{div}(\Phi)(\phi_{-s})\,ds\right),$$

*and for all $p > 1$ with the conjugate exponent $q$, that is, $\frac{1}{p} + \frac{1}{q} = 1$, we have*

$$\|K_t\|_{L^p}^p \leq E[e^{pqt|\mathrm{div}(\Phi)|}].$$

PROOF.   See [4] and [24].   □

Now we are ready to prove the Cameron–Martin theorem on the loop space.

THEOREM 7.3.   *For any $h \in H_0$, the pinned Wiener measure $\nu$ on $L_{m_0}(M)$ is quasi-invariant under Driver's flow $\widetilde{\Phi}_t$. Let*

$$K_t = \frac{d(\widetilde{\Phi}_t)_* \nu}{d\nu}.$$

*Then*

$$(7.7) \quad K_t(\gamma) = \exp\left(\int_0^t \delta(h)(\widetilde{\Phi}_{-s}(\gamma))\, ds\right), \qquad \nu\text{-}a.s.\ \gamma \in L_{m_0}(M),$$

*and for all $p > 1$ with the conjugate exponent $q$, that is, $\frac{1}{p} + \frac{1}{q} = 1$, we have*

$$(7.8) \qquad \qquad \|K_t\|_{L^p(\nu)}^p \leq E_\nu[e^{pqt|\widetilde{\delta}(h)|}].$$

PROOF.   Set $\Omega = L_{m_0}(M)$, $P = \nu$, $\mathcal{C} = \mathcal{F}C(L_{m_0}(M))$ [the collection of all cylindrical functionals on $L_{m_0}(M)$], and let $\phi_t$ be Driver's flow $\widetilde{\Phi}_t$. Theorem 4.2 shows that the divergence $\mathrm{div}(\Phi)$ in (7.5) associated to the flow $\phi_t$ is just $\delta(h)$. By Theorem 5.4, the $\nu$-exponential integrability (7.6) holds for $\delta(h)$. Hence Cruzeiro's lemma applies to $(L_{m_0}(M), \widetilde{\Phi}_t, \nu)$. Thus, $(\widetilde{\Phi}_t)_* \nu$ is absolutely continuous with respect to $\nu$, that is, the pinned Wiener measure $\nu$ is quasi-invariant on loop space $L_{m_0}(M)$ under the flow $\widetilde{\Phi}_t$. Moreover, we obtain (7.7) and the $L^p$-inequality (7.8).   □

According to Remark 1.1, our main results apply to the special case where $M = G$ is a compact connected Lie group equipped with an *Ad*-invariant metric and the left or the right Cartan connection. Theorem 1.5 recaptures the well-known result due to Malliavin and Malliavin [25] on the quasi-invariance of the pinned Wiener measure on the loop group. Indeed, our method is inspired by [25] where the authors initiated the so-called *localization method from paths to loops* based on quasi-sure analysis.

**8. Stochastic anti-development of pinned Brownian motions.**   Since [2], it has been well known that the stochastic anti-development of the pinned Brownian motion on any compact Riemannian manifold is a semimartingale up to time 1. However, the $L^p(\nu)$-convergence $(p \geq 1)$ and the Fernique type exponential integrability theorem for the stochastic anti-development of pinned Brownian motions were first proved by Gross on a compact Lie group (see Lemma 4.8, Remark 4.9 and Corollary 4.10 in [13]). More precisely, let $L_e(G)$ be the loop group over a compact Lie group $G$ equipped with an $Ad$-invariant metric and the left Cartan (or the right Cartan) connection, with $e$ its unit element. The anti-development of the pinned Brownian motion $\{g(s), s \in [0,1]\}$ is given by

$$b(s) = \int_0^s g^{-1}(r) \circ dg(r), \qquad s \in [0,1].$$

The Gross theorem says that there exists a small $\lambda_0$ such that, for all $\lambda < \lambda_0$, we have $E_\nu[\exp(\lambda \max_{s \in [0,1]} \|b(s)\|^2)] < +\infty$. See also [12] for an alternative proof. In this section we will use Lemma 5.2 to establish the $L^p(\nu)$-convergence and a Fernique-type exponential integrability theorem for the stochastic anti-development of pinned Brownian motion on any compact Riemannian manifold with a TSS connection. Our result also sharpens the exponential exponent $\lambda_0$ appearing in the Gross–Fernique theorem. We begin with the following theorem:

THEOREM 8.1. *Let $G$ be a compact Lie group equipped with an Ad-invariant metric and the left or the right Cartan connection. Then for all*

$$\lambda < \lambda_0 = \tfrac{1}{2} \inf\{\|w\|_H^2 : w \in X, \|w\|_{2m,\alpha} = 1\},$$

*we have*

$$E_\nu[\exp(\lambda \|b\|_{2m,\alpha}^2)] < +\infty.$$

*Moreover, for any $\lambda < \tfrac{1}{2}$, we have*

$$E_\nu\left[\exp\left(\lambda \max_{s \in [0,1]} \|b(s)\|^2\right)\right] < +\infty.$$

PROOF. For $\mu$-a.s. $x \in X$, since $dg_x(s) = g_x(s) \circ dx(s)$, we have

$$b(s)(x) = \int_0^1 [g_x(r)]^{-1} \circ dg_x(r) = x(s) \in W^{r,p}(X, T_e G), \qquad s \in [0,1].$$

By the Donsker–Varadhan [5] refinement version of the well-known Fernique lemma, we have

$$E_\mu[\exp(\lambda \|b\|_{2m,\alpha}^2)] = E_\mu[\exp(\lambda \|x\|_{2m,\alpha}^2)] < +\infty,$$

provided that

$$\lambda < \lambda_0 = \inf\{I(w) : w \in X, \ \|w\|_{2m,\alpha} = 1\},$$

where

$$I(w) = \tfrac{1}{2} \|w\|_H^2 \mathbb{1}_{[w \in H]} + \infty \mathbb{1}_{[w \notin H]}.$$

Note that the function $\|\cdot\|_{2m,\alpha}^2$ is smooth in $X^{2m,\alpha}$ in the sense of Fréchet–Gâteaux. Thus, the Wiener functional $x \to \|x\|_{2m,\alpha}^2$ belongs to $W^{\infty,\infty}(X^{2m,\alpha})$. Hence Lemma 5.2 applies to $F(x) = \|x\|_{2m,\alpha}$. That is to say, for any $\varepsilon > 0$ and any $p > 1$, we have

$$E_\nu[\exp(\lambda \|b\|_{2m,\alpha}^2)]$$
$$\leq \{E_\mu[\exp((1+\varepsilon)p\lambda\|b\|_{2m,\alpha}^2)]\}^{1/(p(1+\varepsilon))} \|\|x\|_{2m,\alpha}\|_{W^{n,n^2p\varepsilon/(1+\varepsilon)}(X)}^n,$$

where $n = 2[\frac{l}{2}] + 2$ if we assume that $G \subset \mathbb{R}^l$ is a Nash–Whitney embedding. Thus, for all $\lambda < \frac{\lambda_0}{(1+\varepsilon)p}$, we have $E_\nu[\exp(\lambda\|b\|_{2m,\alpha}^2)] < +\infty$. Since $\varepsilon > 0$ and $p > 1$ are arbitrary, for all $\lambda < \lambda_0$ we get $E_\nu[\exp(\lambda\|b\|_{2m,\alpha}^2)] < +\infty$. Now $\|w\|_\infty := \max_{s \in [0,1]} \|w(s)\|_\infty \leq C\|w\|_{2m,\alpha}$. Thus

$$E_\nu\left[\exp\left(\lambda \max_{s \in [0,1]} \|b(s)\|^2\right)\right] < +\infty$$

holds provided that

$$\lambda \sup_{b \neq 0}\left[\frac{\|b\|_\infty^2}{\|b\|_{2m,\alpha}^2}\right] \leq \lambda \sup_{w \in X \setminus \{0\}}\left[\frac{\|w\|_\infty^2}{\|w\|_{2m,\alpha}^2}\right] < \lambda_0.$$

Set $X^* = X \setminus \{0\}$. Then, for all

$$\lambda < \frac{1}{2} \inf_{w \in X^*}\left[\frac{\|w\|_H^2}{\|w\|_{2m,\alpha}^2}\right] \inf_{w \in X^*}\left[\frac{\|w\|_{2m,\alpha}^2}{\|w\|_\infty^2}\right]$$
$$\leq \frac{1}{2} \inf_{w \in X^*}\left[\frac{\|w\|_H^2}{\|w\|_{2m,\alpha}} \cdot \frac{\|w\|_{2m,\alpha}^2}{\|w\|_\infty^2}\right]$$
$$= \frac{1}{2} \inf_{w \in X^*}\left[\frac{\|w\|_H^2}{\|w\|_\infty^2}\right] = \frac{1}{2},$$

that is, for all $\lambda < \frac{1}{2}$, we have

$$E_\nu\left[\exp\left(\lambda \max_{s \in [0,1]} \|x(s)\|^2\right)\right] < +\infty.$$

The proof of Theorem 8.1 is complete. □

Similarly to the proof of Theorem 8.1, if we replace $b(s) = \int_0^s g^{-1}(r) \circ dg(r)$ by

$$x(s) = \int_0^s U_r^{-1}(\gamma) \circ d\gamma(r), \qquad \nu\text{-a.s. } \gamma \in L_{m_0}(M),$$

then we can prove the following Fernique-type exponential integrability theorem for the stochastic anti-development of pinned Brownian motions on a compact Riemannian manifold.

THEOREM 8.2. *Let $M$ be a compact Riemannian manifold equipped with a torsion-skew symmetric (TSS) connection. Then for all*

$$\lambda < \lambda_0 := \tfrac{1}{2}\inf\{\|w\|_H^2 : w \in X, \ \|w\|_{2m,\alpha} = 1\},$$

*we have*

$$E_\nu\left[\exp\left(\lambda \|x\|_{2m,\alpha}^2\right)\right] < +\infty.$$

*Moreover, for all $\lambda < \tfrac{1}{2}$, we have*

$$E_\nu\left[\exp\left(\lambda \max_{s\in[0,1]} \|x(s)\|^2\right)\right] < +\infty.$$

Finally, let us prove the $L^p(\nu)$-convergence of the stochastic anti-development of pinned Brownian motion.

THEOREM 8.3. *Let $M$ be a compact Riemannian manifold equipped with a TSS connection. Then for any $p \geq 1$, $x(s)$ converges to $x(1)$ in $L^p(\nu)$ as $s$ tends to 1. Moreover, for any $p \geq 1$, there exists a constant $C_p$ such that*

$$\|x(s) - x(1)\|_{L^p(\nu)} \leq C_p (1-s)^{1/2}.$$

PROOF. Similarly to the proof of Theorem 4.1, there exist a constant $C$ and a pair $(r,q) \in \mathbb{N} \times (1,+\infty)$ such that, for $p = 2n$, $n \in \mathbb{N}$,

$$\|x(s) - x(1)\|_{L^p(\nu)}^p \leq C \|\|x(s) - x(1)\|^p\|_{W^{r,q}(X)}.$$

With respect to the Wiener measure on the Wiener space $X$, $x(s) - x(1)$ is a centered Gaussian variable with variance $1 - s$. Thus, there exists a constant $C(n,q)$ such that

$$E_{\mu_0}[\|x(s) - x(1)\|^{2qn}] \leq C(n,q)(1-s)^{nq}.$$

On the other hand, for all $h_1,\ldots,h_r \in H$, $i = 1,\ldots,d$, we have

$$\nabla_{h_1} \cdots \nabla_{h_r}(x_i(s) - x_i(1))^{2n} = 2n(2n-1)\cdots(2n-r)(x_i(s) - x_i(1))^{2n-r}$$
$$\times (h_1^i(s) - h_1^i(1))\cdots(h_r^i(s) - h_r^i(1)),$$

from which one can easily verify that

$$\|\nabla^r \|x(s) - x(1)\|^{2n}\|_q \leq C(n,r,q)(1-s)^{(2n-r)/2}.$$

Therefore, we get

$$\|x(s) - x(1)\|_{L^{2n}(\nu)} \leq C(n,r,q)[(1-s)^{1/2} + (1-s)^{(2n-r)/4n}].$$

Note that $r$ and $q$ are independent of $n$. Hence

$$\|x(s) - x(1)\|_{L^{2n}(\nu)} \leq C_n (1-s)^{1/2}.$$

This yields that $\|x(s) - x(1)\|_{L^p(\nu)} \leq C_p (1-s)^{1/2}$ for all $p \geq 1$. □

**9. Two remarks on Doob's $h$-processes approach.** For the completeness of the paper, we would like to give two remarks on the Doob $h$-processes approach for studying the problems discussed in this paper.

Recall that with respect to the pinned Wiener measure $\nu$ on loop space $L_{m_0}(M)$, the conditional Brownian motion $\gamma_t$ is given by the following Stratonovich SDE:

$$d\gamma_s = U_s \circ db_s, \qquad \gamma_0 = m_0,$$

where $\{b_s, s \in [0,1]\}$ is the anti-development of $\{\gamma_s, s \in [0,1]\}$ through the Itô map $I : X \to P_{m_0}(M)$ which can be well defined up to a slim subset of $X$. Moreover, $\{b_s, s \in [0,1]\}$ is a semimartingale with the following Doob-Meyer decomposition:

$$db_s = d\beta_s + U_s^{-1} \nabla \log p_{1-s}(\gamma_s, m_0)\, ds,$$

where $\{\beta_s, s \in [0,1]\}$ is a $\nu$-Brownian motion on $(X, \mathcal{F}_s, \mathcal{F}, \nu)$. See, for example, [2, 7, 9, 15, 16].

Let $h \in H$ and let $D_h$ be the vector field on $L_{m_0}(M)$. Let $\gamma^t$ be the Driver flow on the loop space $L_{m_0}(M)$ given by

$$\dot{\gamma}^t = D_h(\gamma^t) = U_s^t h(s),$$
$$\gamma^0 = \gamma,$$

where $\{U_s^t, s \in [0,1]\}$ is the horizontal lift of $\{\gamma_s^t, s \in [0,1]\}$. Using the intertwining formula for the differential of the stochastic development map $I^{-1} : P_{m_0}(M) \to X$ (see [6, 8, 10, 15, 22, 24]), the pull-back flow $b^t = I^{-1}(\gamma^t)$ satisfies the following equation:

$$\frac{\partial}{\partial t} d_s b_s^t = \dot{h}(s)\, ds - q_h(\gamma^t, s) \circ d_s b_s^t,$$
$$b_s^0 = b_s,$$

where

$$d_s b_s^t = d_s \beta_s^t + [U_s^t]^{-1} \nabla \log p_{1-s}(\gamma_s^t, m_0)\, ds.$$

Thus, we have

$$\left[\frac{\partial}{\partial t} + q_h(\gamma^t, s)\right] d_s b_s^t = \dot{h}_s ds,$$

whence

$$\left[\frac{\partial}{\partial t} + q_h(\gamma^t, s)\right] \circ d_s \beta_s^t$$
$$= \dot{h}_s ds - \left[\frac{\partial}{\partial t} + q_h(\gamma^t, s)\right]([U_s^t]^{-1} \nabla \log p_{1-s}(\gamma_s^t, m_0))\, ds$$

$$= \dot{h}_s ds - q_h(\gamma^t, s)([U_s^t]^{-1} \nabla \log p_{1-s}(\gamma_s^t, m_0)) \, ds$$
$$+ [U_s^t]^{-1} \left(\frac{\partial}{\partial t} U_s^t\right) [U_s^t]^{-1} \nabla \log p_{1-s}(\gamma_s^t, m_0) \, ds$$
$$- [U_s^t]^{-1} \nabla^2 \log p_{1-s}(\gamma_s^t, m_0) \left(\frac{\partial}{\partial t} \gamma_s^t\right) ds.$$

Moreover, by the Bismut formula (see, e.g., [2, 22]), we have

$$\frac{\partial}{\partial t} U_s^t = U_s^t q_h(\gamma^t, s).$$

Combining the above formulas, we can derive the pull-back Driver flow equation as follows:

$$\left[\frac{\partial}{\partial t} + q_h(\gamma^t, s)\right] \circ d_s \beta_s^t = \dot{h}_s ds - [U_s^t]^{-1} \nabla^2 \log p_{1-s}(\gamma_s^t, m_0) U_s^t h(s) \, ds.$$

In order to use the standard Picard iteration or the Euler iteration method to solve the above flow equation, it is clear that one has to use the Hessian estimate of the logarithm of the heat kernel on compact Riemannian manifold. Moreover, in order to use the usual approach based on the Girsanov theorem and Lévy's invariance of Brownian motion under adapted rotations to prove the quasi-invariance of the pinned Wiener measure under the pull-back Driver flow, we need to verify the Novikov exponential integrability condition of the drift term given in the right-hand side of the flow equation for $\beta^t$. Thus, we need to use again the Hessian estimate of the logarithm of the heat kernel. See, for example, [7, 9, 15, 16]. See also [17, 27, 29, 31] for the short time estimates of logarithmic derivatives of the heat kernel.

To end this paper, let us mention that Gong has informed us that, by using the gradient and the Hessian estimates of the logarithm of the heat kernel, Gong and Ma can also prove the $L^p(\nu)$-convergence and the $\nu$-exponential integrability of the divergence functional $\delta(h)$ (for all $\lambda < \lambda_0$ for some constant $\lambda_0$ which depends on $\|h\|_H$ and possibly on the constants appearing in the gradient and the Hessian estimates of the logarithmic heat kernel) by the Doob $h$-processes approach. Without using heat kernel estimate, our approach based on the Airault–Malliavin–Sugita–Watanabe inequality (see Theorem 2.2 and Lemma 5.2) and Sobolev estimates shows that we can get an explicit estimate for $\lambda_0$ which only depends on the uniform bound of the Ricci curvature $\|\text{Ric}\|_\infty$ and $\|h\|_H$, that is, $\lambda_0 = [(2 + \|\text{Ric}\|_\infty)\|h\|_H]^{-1}$.

**Acknowledgments.** The earliest version of this paper was done in 1996 and appeared in the author's Ph.D. thesis [21] under the supervision of Professors Zhi-Ming Ma, P. Malliavin and A. B. Cruzeiro. The author thanks his supervisors and Professor Shizan Fang for suggesting to him the use of quasi-sure analysis to study the problems in this paper. Thanks are also due to Professors B. Driver, F. Z. Gong, E. Hsu, M. Ledoux, R. Léandre, T. Lyons, J. Ren and D. W. Stroock for helpful discussions and constant encouragements. Finally, the author would like to thank the Associate Editor and the anonymous referee for their careful reading, comments for correction and helpful suggestions for improvements.

MATHEMATICAL INSTITUTE
UNIVERSITY OF OXFORD
24-29 ST. GILES'
OXFORD, OX1 3LB
UNITED KINGDOM
AND
LABORATOIRE DE STATISTIQUE ET PROBABILITÉS
UNIVERSITÉ PAUL SABATIER
118 ROUTE DE NARBONNE
31062 TOULOUSE CEDEX
FRANCE
E-MAIL: xiang@cict.fr